\documentclass[11pt,reqno,tbtags]{amsart} %{article} %
\usepackage[utf8]{inputenc}

\usepackage[a4paper,width=150mm,top=25mm,bottom=25mm]{geometry}
\usepackage{mathtools}
\usepackage{suffix}
\usepackage{enumerate}
\usepackage{enumitem}
\usepackage{listings}

\makeatletter
\newcommand*{\rom}[1]{\expandafter\@slowromancap\romannumeral #1@}
\makeatother
\usepackage{pgf, tikz}
\usepackage{subcaption}
\usetikzlibrary{arrows, automata}
\usepackage{float}

\usepackage{parskip}
\usepackage{amsmath,amsthm,amssymb}
\numberwithin{equation}{section}

\allowdisplaybreaks
\usepackage{bbm}
\usepackage[makeroom]{cancel}

\usepackage{xcolor}
\definecolor{coolblack}{rgb}{0.0, 0.18, 0.39}

\usepackage[breaklinks=true]{hyperref}
\hypersetup{
colorlinks=true,
linkcolor=blue,
urlcolor=blue,
citecolor=black
}

\theoremstyle{plain}
\newtheorem{theorem}{Theorem}[section]

\newtheorem{lemma}[theorem]{Lemma}

\newtheorem{conjecture}[theorem]{Conjecture}

\theoremstyle{definition}

\newtheorem{definition}[theorem]{Definition}
\newtheorem{example}[theorem]{Example}

\newtheorem{remark}[theorem]{Remark}

\newcommand{\al}{\alpha}
\newcommand{\IP}{\mathbbm{P}}
\newcommand{\E}{\mathbb{E}}

\newcommand{\Pa}{\pi_\alpha} 
\newcommand{\Pb}{\pi_\beta}

\newcommand\numberthis{\addtocounter{equation}{1}\tag{\theequation}}

\newcommand{\tone}{\mathbf{1}}

\newcommand{\IZ}{\mathbbm{Z}}

\newcommand{\cE}{\mathcal{E}}

\newcommand{\cT}{\mathcal{T}}
\newcommand{\var}{\mathrm{Var}}
\newcommand{\cov}{\mathrm{Cov}}
\newcommand{\dtv}{\mathop{d_{\mathrm{TV}}}}
\newcommand{\dw}{\mathop{d_{\mathrm{W}}}}
\newcommand{\dk}{\mathop{d_{\mathrm{K}}}}
\newcommand{\law}{\mathcal{L}}

\renewcommand{\leq}{\leqslant}
\renewcommand{\geq}{\geqslant}
\renewcommand{\phi}{\varphi}

\renewcommand{\le}{\leq}
\renewcommand{\ge}{\geq}
\newcommand{\refT}[1]{Theorem~\ref{#1}}
\newcommand{\refTs}[1]{Theorems~\ref{#1}}

\newcommand{\refL}[1]{Lemma~\ref{#1}}

\newcommand{\refR}[1]{Remark~\ref{#1}}

\newcommand{\refS}[1]{Section~\ref{#1}}

\newcommand{\refE}[1]{Example~\ref{#1}}

\newcommand{\refConj}[1]{Conjecture~\ref{#1}}
\newcommand\ga{\alpha}
\newcommand\gb{\beta}
\newcommand\gG{\Gamma}
\newcommand\gL{\Lambda}

\newcommand\cI{\mathcal I}
\newcommand\suma{\sum_{\nu\in\cI}}
\newcommand\sumin{\sum_{i=1}^n}

\newcommand\sumn{\sum_{n=1}^\infty}

\newcommand\set[1]{\ensuremath{\{#1\}}}

\newcommand\xpar[1]{(#1)}
\newcommand\bigpar[1]{\bigl(#1\bigr)}
\newcommand\Bigpar[1]{\Bigl(#1\Bigr)}
\newcommand\biggpar[1]{\biggl(#1\biggr)}
\newcommand\lrpar[1]{\left(#1\right)}
\newcommand\bigsqpar[1]{\bigl[#1\bigr]}
\newcommand\sqpar[1]{[#1]}

\newcommand\lrsqpar[1]{\left[#1\right]}

\newcommand\bigabs[1]{\bigl\lvert#1\bigr\rvert}

\newcommand\upto{\nearrow}
\newcommand{\tend}{\longrightarrow}

\newcommand\pto{\overset{\mathrm{p}}{\tend}}
\newcommand\asto{\overset{\mathrm{a.s.}}{\tend}}
\newcommand\ntoo{\ensuremath{{n\to\infty}}}
\newcommand\Ntoo{\ensuremath{{N\to\infty}}}

\newcommand\Be{\operatorname{Be}}

\renewcommand\P{\IP}
\newcommand\Var{\operatorname{Var}}
\newcommand\Cov{\operatorname{Cov}}
\newcommand\xdots{\cdots}
\newcommand\xnot{\text{not }}
\newcommand\bbN{\mathbb N}
\newcommand\jq{q}
\newcommand\gab{\ga\gb}
\newcommand\gaxb{\ga{\cdot}\gb}
\newcommand\gabxcc{{\ga{\cdot}\gb{*}\gamma_1\gamma_2}}
\newcommand\gabcc{{\ga{\cdot}\gb{\cdot}\gamma_1\gamma_2}}
\newcommand\gaxcc{{\ga{*}\gamma_1\gamma_2}}
\newcommand\gacc{{\ga{\cdot}\gamma_1\gamma_2}}
\newcommand\ttone{\tilde{\tone}}
\newcommand\gam{\gamma}

\newcommand\gss{\sigma^2}
\newcommand\aaa{^{(a)}}
\newcommand\logn[1]{\log^{#1}n}

\newcommand\WWn[1]{W_n^{(#1)}}
\newcommand\WIJ[1]{W'_{j,i,#1}}
\newcommand\WKJ[1]{W'_{j,k,#1}}
\newcommand\Sl{S_\ell}
\newcommand\Sli{S_{\ell-1}}

\newcommand\ddx{\mathrm{d}}
\DeclarePairedDelimiter\ceil{\lceil}{\rceil}
\DeclarePairedDelimiter\floor{\lfloor}{\rfloor}

\def\given{\typeout{Command 'given' should only be used within bracket command}}
\newcounter{@bracketlevel}
\def\@bracketfactory#1#2#3#4#5#6{
\expandafter\def\csname#1\endcsname##1{%
\addtocounter{@bracketlevel}{1}%
\global\expandafter\let\csname @middummy\alph{@bracketlevel}\endcsname\given%
\global\def\given{\mskip#5\csname#4\endcsname\vert\mskip#6}\csname#4l\endcsname#2##1\csname#4r\endcsname#3%
\global\expandafter\let\expandafter\given\csname @middummy\alph{@bracketlevel}\endcsname
\addtocounter{@bracketlevel}{-1}}%
}
\def\bracketfactory#1#2#3{%
\@bracketfactory{#1}{#2}{#3}{relax}{1mu plus 0.25mu minus 0.25mu}{0.6mu plus 0.15mu minus 0.15mu}
\@bracketfactory{b#1}{#2}{#3}{big}{1mu plus 0.25mu minus 0.25mu}{0.6mu plus 0.15mu minus 0.15mu}
\@bracketfactory{bb#1}{#2}{#3}{Big}{2.4mu plus 0.8mu minus 0.8mu}{1.8mu plus 0.6mu minus 0.6mu}
\@bracketfactory{bbb#1}{#2}{#3}{bigg}{3.2mu plus 1mu minus 1mu}{2.4mu plus 0.75mu minus 0.75mu}
\@bracketfactory{bbbb#1}{#2}{#3}{Bigg}{4mu plus 1mu minus 1mu}{3mu plus 0.75mu minus 0.75mu}
}
\bracketfactory{clc}{\lbrace}{\rbrace}
\bracketfactory{clr}{(}{)}
\bracketfactory{cls}{[}{]}
\bracketfactory{abs}{\lvert}{\rvert}
\bracketfactory{norm}{\Vert}{\Vert}
\bracketfactory{floor}{\lfloor}{\rfloor}
\bracketfactory{ceil}{\lceil}{\rceil}
\bracketfactory{angle}{\langle}{\rangle}

\usepackage[maxbibnames=99,style=numeric,backend=biber,sorting=nyt, bibstyle=authoryear]{biblatex}

\DeclareFieldFormat
  [article,book,incollection,inproceedings,misc,thesis,unpublished]
  {title}{#1\isdot}

\AtNextBibliography{\small}
\DeclareNameAlias{sortname}{family-given}
\DeclareBibliographyCategory{}

\DeclareFieldFormat{labelnumberwidth}{\mkbibbrackets{#1}}
\defbibenvironment{bibliography}
  {\list
     {\printtext[labelnumberwidth]{%
    \printfield{prefixnumber}%
    \printfield{labelnumber}}}
     {\setlength{\labelwidth}{\labelnumberwidth}%
      \setlength{\leftmargin}{\labelwidth}%
      \setlength{\labelsep}{\biblabelsep}%
      \addtolength{\leftmargin}{\labelsep}%
      \setlength{\itemsep}{\bibitemsep}%
      \setlength{\parsep}{\bibparsep}}%
      }
  {\endlist}
  {\item}

\renewbibmacro{in:}{%
  \ifentrytype{article}
    {}
    {\bibstring{in}%
     \printunit{\intitlepunct}}}
     
\usepackage{xpatch}
\xpatchbibmacro{cite}{\printnames{labelname}}%
{\ifciteseen{\printnames{labelname}}{\printnames[][-\value{listtotal}]{labelname}}}
{}{}

\xpatchbibmacro{textcite}{\printnames{labelname}}%
{\ifciteseen{\printnames{labelname}}{\printnames[][-\value{listtotal}]{labelname}}}
{}{}
\makeatletter
\patchcmd{\@maketitle}{\LARGE \@title}{\fontsize{16}{19.2}\selectfont\@title}{}{}
\makeatother

\title{Approximation of subgraph counts in the uniform attachment model}
\author{Johan Bj\"orklund, Cecilia Holmgren, Svante Janson, Tiffany Y.\ Y.\ Lo}
\thanks{Supported by the Knut and Alice Wallenberg Foundation, %SJ, CH, TL?
 Ragnar Söderberg's Foundation, %CH
the Swedish Research Council,  %CH, TL?
and the David Lachlan Hay Memorial Fund%  %TL
}
\date{7 November 2023} %; compiled \today}

\begin{document}
\begin{abstract}
 We use Stein's method to obtain distributional approximations of
    subgraph counts in the uniform attachment model or random directed
    acyclic graph; we provide also estimates of rates of convergence. 
In particular,
    we give uni- and multi-variate Poisson approximations to the 
    counts of cycles, and normal approximations to the counts of unicyclic
    subgraphs; we also give a partial result for the counts of trees. 
We further find a class of multicyclic graphs
 whose subgraph counts are a.s.\ bounded as $n\to\infty$. 
\end{abstract}

\maketitle
\section{Introduction}
A random uniform attachment graph, denoted $G^m_n$ and also known as a 
uniform random recursive \emph{dag} (directed acyclic graph), 
can be constructed recursively as follows. Fix
$m\geq 1$.
The initial graph $G^m_1$ is a single isolated
vertex. To construct $G^m_{n}$ from $G^m_{n-1}$, we add vertex $n$ to
$G^m_{n-1}$, with vertex $n$ born with $m$ edges, 
which are labelled by $1,\dots,m$. 
For %$1\leq i\leq m$, 
$i\in[m]:=\set{1,\dots,m}$,
the
other endpoint of $n^{(i)}$, the $i$th edge of vertex $n$, is then uniformly
chosen among the existing vertices of $G^m_{n-1}$,
i.e., in $[n-1]$. Thus $G^m_n$ has $n$ vertices and $(n-1)m$ edges, and each edge has a label
in $[m]$.
Observe that we allow multiple edges when $m\geq 2$. 
An edge in $G^m_n$ can be thought of as
always pointing towards the vertex with the smaller label, and so there is
no real distinction between the undirected and directed versions of the uniform
attachment graph.  

When $m=1$, the model is the random recursive tree first studied in \textcite{na1970distribution}; and when $m\geq 2$, the model was first introduced in \textcite{devroye1995strong}. There is an abundance of literature on random recursive trees (see e.g.\ \textcite{drmota2009} for an overview), but here we mention \textcite{holmgren2015limit, feng2010, feng2008phase, fuchs2008}, which provide Poisson and normal approximations to the counts of subtree copies. For $m\geq 2$, results on vertex degrees can be found in \textcite{devroye1995strong,tsukiji2001}, and results on depths and path lengths are available in \cite{arya1999,broutin2012,devroye2011,tsukiji1996}. The recent paper \textcite{janson2023} studies the number of vertices that can be reached from vertex $n$ via a directed path, where the edge is thought as pointing from the larger vertex label towards the smaller one.

In this article, we consider the counts of fixed subgraphs (not depending on
$n$)  as $n\to\infty$;
the parameter $m$ is assumed to be fixed with $m\ge2$. 
We provide uni- and multivariate
Poisson approximations to the counts of cycles. We also prove that the
counts of unicyclic subgraphs are asymptotically normal (after suitable
renormalisation). We conjecture that the counts of more general trees are
approximately normal, and show that this conjecture holds for 
stars.
All these approximation results are accompanied with convergence
rates. We also consider multicyclic subgraphs. In particular, we show that
if in addition to having more than two cycles, the subgraph is `leaf-free'
(i.e.\ has no vertices of degree one), then the number of copies in the
uniform attachment graph is bounded a.s. For multicyclic subgraphs that
are not leaf-free, we identify the exact rate of growth of the expectation
and an upper bound on the variance. We conjecture that the count of a
subgraph of this type converges a.s.\ after suitable rescaling, and prove a
special case.

The distributional approximation tool that we use here is Stein's method, using in particular the size-bias coupling for both types of approximation; see \cite{barbour1992poisson, chen2010normal, ross2011fundamentals,} for an overview. This is in contrast to the analytic and contraction methods employed in \textcite{feng2010, feng2008phase, fuchs2008}. In  
\cite{holmgren2015limit}, Stein's method was applied in proving the Poisson and normal approximation results for the counts of subgraphs and their functionals in random recursive trees. However, we note that, for normal approximation, the approach that we use here is different from that in \cite{holmgren2015limit}, which does not provide a convergence rate.

A major challenge in the distributional approximation problem of subgraph
counts in the uniform attachment graph is the computation of the 
variance of the counts. In particular, our methods require a lower bound
on the variance.
In the case $m=1$, the variance can be computed
explicitly by using an elegant bijection between the random binary tree and
the random recursive tree \cite[Section 2.3.2]{knuth}, as done in
\cite{holmgren2015limit}. When $m\geq 2$, the order of the variance of the
count of certain non-tree subgraphs  
can be obtained by analysing the 
covariances of the indicators that a certain subgraph exists in $G^m_n$,
and finding the
pairs of subgraphs that contribute the
dominant term in the variance, but this task becomes difficult in the case
of  trees. In contrast to the case of $m=1$, it is also much harder to
obtain an explicit expression for the mean and variances in the case of
$m\geq 2$. This is due to the fact for most subgraphs, the same set of
vertices in the uniform attachment graph can form a copy in several
ways. Consequently, we are only able to give the order of
the mean and variance (i.e., within constant factors),
and in some cases only an upper bound, 
and we only consider fixed subgraphs and a fixed $m$.

\begin{remark}
It would be interesting to extend the results of this paper to subgraphs
$H_n$ depending on $n$, or to $m=m(n)$ growing with $n$.
Our methods still work in principle, but it becomes more complicated to
estimate the various expressions, and we have not pursued this. 
We leave this as open problems.
\end{remark}

\section{Notation}\label{Snot}

The objective of this paper is thus to study the number of copies of a given
graph $H$ in $G^m_n$. 
We will only consider connected $H$.
We emphasise that we generally see both $G^m_n$ and
the $H$ as \emph{undirected}. 
However, each edge in $G^m_n$ has a canonical direction where an edge
between the vertices $i$ and $j$ is directed towards the smaller of $i$ and
$j$. Thus every copy of $H$ in $G^m_n$ also has an induced direction on its
edges. Note, however, that different copies may induce different directions in
$H$.

Let $K^m_n$
be the \emph{$m$-fold complete multigraph on $n$},
defined 
to be the multigraph on $[n]$ such that
for every pair $(i,j)$ with $1\le i< j\le n$, there are
$m$ edges $j\to i$, and these edges are labelled $1,\dots,m$.
We denote the edge $j\to i$ with label $a$ by $ji^{(a)}$.
Note that we can regard $G^m_n$ as a subgraph of $K^m_n$ in the obvious way
(preserving edge labels);
in fact, we may define $G^m_n$ as the subgraph of $K^m_n$ obtained by
choosing, for every $j\in\{2,\dots,n\}$ and every $a\in[m]$,
exactly one of the edges $ji\aaa$ with label $a$ and $i\in[j-1]$,
all choices uniformly at random and independent.

Let $\gG$ be the set of all copies of $H$ in $K^m_n$  that 
do not contain two different edges $ji_1^{(a)}$ and $ji_2^{(a)}$ with $i_1,i_2<j$ and the
same label $a$. (If there are two such edges $ji_1^{(a)}$ and $ji_2^{(a)}$,
then the copy can never be a subgraph of $G^m_n$.)
We say that $\gG$ is  the set of \emph{potential copies} of $H$ in $G^m_n$.
For $\ga\in\Gamma$, let $\tone_\al$ be the
indicator 
variable that takes value 1 if $\ga\subseteq G^m_n$,
and 0 otherwise. 
The number $W_n$ of  copies of $H$ in $G^m_n$ is thus
\begin{align}
W_n=
\sum_{\al\in\Gamma}\tone_\al,  
\end{align}
and 
we are therefore interested in approximating the distribution of this sum.

\begin{example}\label{E3}
For a simple example, let $H$ be a triangle.
Then, 
a potential copy of $H$ in $G^m_n$ is described by its vertex set
$\set{i,j,k}\subseteq[n]$, where we may assume 
$i<j<k$, together with its edges $ji^{(a)}$, $ki^{(b)}$, and $kj^{(c)}$,
where $a,b,c\in[m]$ and $b\neq c$. 
\end{example}

\begin{remark}\label{RmH}
  For any graph $H$ there is an integer $m_H\ge1$ such that if 
$m<m_H$, then for every $n$ 
there are no potential copies of $H$ in $G^m_n$, so $\gG=\emptyset$
and thus $W_n=0$ deterministically,
but if $m\ge m_H$ and $n$ is large enough, then there are potential copies.
For example, $m_H=1$ if $H$ is a tree, $m_H=2$ if $H$ is a cycle, and
$m_H=k-1$ if $H$ is the complete graph on $k$ vertices.
\end{remark}

\begin{remark}\label{Rmulti}
Recall that $G^m_n$ is a loop-less multigraph. Similarly, 
$H$ can be a loop-less multigraph in the results below,
although we for simplicity  write ``graph'' and ``subgraph''.
\end{remark}

\subsection{Probability distances}
To precisely state our results, we also need to define the metrics in
consideration here. The \emph{total variation distance} $\dtv$
between two probability measures $\nu_1$ and $\nu_2$ supported on $\IZ^+$ is defined as 
\begin{align}
    \dtv(\nu_1,\nu_2)=\sup_{A\subset \IZ^+}|\nu_1(A)-\nu_2(A)|;
\end{align}
and the \emph{Wasserstein distance} $\dw$ between two probability measures $\nu_1,\nu_2$ supported on $\mathbbm{R}$ is defined as 
\begin{equation}
    \dw(\mu,\nu)=\sup_f \biggl|\int f(x)d\mu(x)-\int f(x)d\nu(x)\biggr|, 
\end{equation}
where the supremum is over all $f:\mathbbm{R}\to \mathbbm{R}$ such that $|
f(x)-f(y)|\leq | x-y|$ for any $x,y\in \mathbbm{R}$. Note that 
if, for example, $Y$ is a random variable with the standard normal
distribution $\mathcal{N}(0,1)$, 
then for any random variable $X$, the usual
Kolmogorov distance can be bounded by
\begin{align}\label{KW}
d_{\mathrm{K}}\bigpar{\law(X),\law(Y)}:=\sup_{a\in \mathbbm{R}}\bigabs{\IP[X\leq a]-\IP[Y\leq a]}
=
O\bigpar{\sqrt{\dw(\law(X),\law(Y)}}
\end{align}
(see e.g.\ \cite[Proposition 2.1]{ross2011fundamentals}).

\subsection{Further notation}
We sometimes tacitly assume that $n$ is not too small.

$\asto$ and $\pto$ denote convergence almost surely (a.s.)\
and in probability, respectively.

$\mathrm{Po}(\mu)$ denotes the Poisson distribution with mean $\mu$,
and $\mathcal{N}(0,1)$ is the standard normal distribution.

$C$ denotes constants that may vary from one occurrence to the
next. 
They do not depend on $n$, but they may depend on $m$, $H$, and other
parameters.

\section{Main results}
Consider first subgraphs that are cycles.
The first two theorems show that when $n$ is large, the count of any 
cycle with a fixed number of edges 
is approximately Poisson, 
and that the joint
distribution of the counts of cycles of different numbers of edges
can be approximated by independent Poisson variables. For convenience, we
refer to
a cycle with $\ell$ edges and vertices simply as an $\ell$-cycle. Note also
we view a pair of parallel edges as a 2-cycle. 
\begin{theorem}\label{th:pocycle}
Fixing the positive integers $m\ge2$ and $\ell\geq 2$, let $W_n$ be the
number of $\ell$-cycles in $G^m_n$, and let $\mu_n:=\E W_n$. Then, 
\begin{gather}\label{eq:cyclem}
    \mu_n =\Theta(\log n)
 \end{gather}
and there is a positive constant $C=C(m,\ell)$ such that
\begin{align}\label{th:unipo}
    \dtv(\law(W_n),\mathrm{Po}(\mu_n))\leq C\log^{-1} n.
\end{align}
\end{theorem}

\begin{remark}\label{re:PoNaCycle}
    Let $W_n$ be as above and $Y'_n\sim \mathrm{Po}(\mu_n)$, and define
    \begin{equation}\label{cycountsc}
        Z_n:=\frac{W_n-\mu_n}{\sqrt{\mu_n}},\qquad Y_n:=\frac{Y'_n-\mu_n}{\sqrt{\mu_n}}.
    \end{equation}
    It follows from \refT{th:pocycle} that
    \begin{equation}\label{kolcycle}
        d_{\mathrm{K}}\bigpar{\law(Z_n), \law(Y_n)}\leq \dtv\bigpar{\law(Z_n), \law(Y_n)} \leq C\logn{-1}.
    \end{equation}
    On the other hand, the classical Berry--Esseen theorem implies that 
    \begin{equation}\label{be}
        d_{\mathrm{K}}\bigpar{\law(Y_n), \mathcal{N}(0,1)}\leq C\mu_n^{-\frac{1}{2}},
    \end{equation}
    Combining \eqref{kolcycle} and \eqref{be} with the triangle inequality, and using \eqref{eq:cyclem}, we obtain
    \begin{equation}\label{kolbdcyc}
    d_{\mathrm{K}}\bigpar{\law(Z_n), \mathcal{N}(0,1)}\leq C \logn{-\frac{1}{2}}.
    \end{equation}
In particular, the cycle count $W_n$ is asymptotically normal.
In fact, the estimate \eqref{kolbdcyc} is sharp.
Since $\mathcal N(0,1)$ has a continuous distribution function, while
the distribution function of $Z_n$ has a jump $\P(W_n=k)$ (if this is
non-zero) at $(k-\mu_n)/\sqrt{\mu_n}$, 
it follows by choosing $k=\floor{\mu_n}$ 
and using \eqref{th:unipo} and \eqref{eq:cyclem} that
\begin{align}\label{bea}
d_{\mathrm{K}}\bigpar{\law(Z_n), \mathcal{N}(0,1)}
%\ge\tfrac12\sup_{k\ge0}\P\bigpar{W_n=k}
\ge\tfrac12\P\bigpar{W_n=\floor{\mu_n}}
=\Theta\bigpar{\mu_n^{-\frac12}}
=\Theta\bigpar{\logn{-\frac{1}{2}}}.  
\end{align}
Consequently, combining \eqref{be} and \eqref{bea},
\begin{align}
d_{\mathrm{K}}\bigpar{\law(Z_n), \mathcal{N}(0,1)}
=\Theta\bigpar{\logn{-\frac{1}{2}}}
.\end{align}
\end{remark}

In the following corresponding theorem on  multi-variate Poisson
approximation, the error bound that we obtain is slightly inferior to the one in
\refT{th:pocycle}.
As in \refR{re:PoNaCycle}, this theorem implies also 
a multivariate normal approximation
(we omit the details).

\begin{theorem}\label{th:mpocycle}
Fix the positive integers $m,r\ge2$ and $\ell(i)\in[2,\infty)$,
$i\in[r]$, 
with $\ell(i)\ne \ell(j)$ for $i\ne j$. Let $W^{(i)}_n$ be the number of
cycles of length $\ell(i)$ in $G^m_n$, and let $\mu_{n,i}:=\E W_n^{(i)}$. Then, there is a positive constant $C=C(m,(\ell(i))_{i\in[r]})$ such that 
\begin{equation}
\dtv\bbclr{\law(\clc{W^{(i)}_n}^r_{i=1}), \prod^r_{i=1}\mathrm{Po}(\mu_{n,i})}\leq \frac{C\log \log n}{\log n}.
\end{equation}
\end{theorem}

\begin{example}
The case $\ell=2$ is simple in \refTs{th:pocycle} and \ref{th:mpocycle}, 
since the numbers of pairs of parallel edges with
larger endpoint $j$ are independent for different $j$;
moreover, we have the exact formula $\mu_n=\binom m2\sum_{j=2}^n\frac{1}{j-1}$. 
For $\ell=3$, \refE{E3} yields
\begin{align}
  \mu_n = m^2(m-1)\sum_{1\le i<j<k\le n}\frac{1}{(j-1)(k-1)^2}
%= m^2(m-1)\sum_{2\le j<k\le n}\frac{1}{(k-1)^2}
%= m^2(m-1)\sum_{3\le k\le n}\frac{k-2}{(k-1)^2}
=m^2(m-1)\bigpar{\log n+ O(1)}.
\end{align}
\end{example}

Next, we state the normal approximation for the count of any
\emph{unicyclic} graph,
i.e., a
graph that contains exactly one cycle.  
For completeness, we
include cycles in the theorem below. Let
$\mathcal{C}_\ell$ be an $\ell$-cycle,
let $s\ge0$, and 
let, for $i=1,\dots,s$, 
$\cT_i$
be a
 tree with $t_i$ edges and a distinguished root, so that $\cT_i$ has
$t_i+1$ vertices. We consider a graph
$\Lambda:=\Lambda_{\ell,t_1,\dots,t_s}$ which can be constructed by
attaching each $\cT_i$ to $\mathcal{C}_\ell$, using a vertex of
$\mathcal{C}_\ell$ as the distinguished root of $\cT_i$. 
(This does not specify $\Lambda$ uniquely, but we choose one possibility.)
By combining the
trees $\cT_i$ that are attached at the same vertex in the cycle, we may
assume that $s\leq \ell$ and that each $\cT_i$ is attached to $\mathcal{C}_\ell$
at a distinct vertex. Denote by $\Gamma=\Gamma^{(n,m)}_{\Lambda}$ the set of
all  potential copies of $\Lambda$; 
%in the uniform attachment graph $G^m_n$;
and saving notation, let 
\begin{align}\label{Wmugs}
  W_n:=\sum_{\al\in\Gamma}\tone_\al,
\qquad
\mu_n:=\E W_n,
\qquad
\sigma_n^2:=\var(W_n).
\end{align}

\begin{theorem}\label{th:nact}
Fixing the integers $m\ge2$, $\ell\ge2$, $s\ge0$, $t_1,\dots,t_s\ge1$,
and a unicyclic subgraph $\Lambda$ as above, let $W_n$, $\mu_n$ and
$\sigma_n$ be as in \eqref{Wmugs},
and $Y_n:=(W_n-\mu_n)/\sigma_n$. Let $t=\sum^s_{i=1}t_i\ge0$; then
\begin{align}\label{eq:ctm}
    \mu_n&=\Theta\bigpar{\logn{t+1}}.
\\ \label{eq:ctv}
\sigma^2_n&=\Theta\bigpar{\logn{2t+1}},
\end{align}
and there is a positive constant $C=C(m,\ell,t)$ such that 
\begin{equation}\label{eq:ctna}
    \dw(\law(Y_n), \mathcal{N}(0,1))\leq C\logn{-\frac{1}{2}}.
\end{equation}
\end{theorem}     

\begin{remark}
    In view of \eqref{KW} and \eqref{eq:ctna}, \refT{th:nact} implies an
    error bound 
    \begin{align}
\dk\bigpar{\law(Y_n), \mathcal{N}(0,1)}=O\bigpar{\logn{-\frac{1}{4}}}      
    \end{align}
for the Kolmogorov distance,
which
in the case of a cycle
     is clearly not as sharp as the error bound in \eqref{kolbdcyc}. 
\end{remark}

The next theorem concerns trees,
and the normal approximation result
relies on the assumption that the variance of the counts of the  tree
of choice is of the exact order $\Theta(n)$. The precise order of the
variance in this case is much harder to establish, essentially due to the
fact that the total covariance of the positively correlated pairs of copies
and that of the negatively correlated pairs of copies are both of order
$O(n)$. However, we are 
able to prove that the variance of the count of a star
is precisely $\Theta(n)$. Below, let $\Lambda$ be a fixed  tree on $t$
vertices,  let $\gG$ be the set of potential copies of $\gL$, 
and  define again $W_n,\mu_n$,
and $\sigma^2_n$ by
\eqref{Wmugs}.

\begin{theorem}\label{th:natree}
Fix the positive integers $m,t$ and the  tree $\Lambda$ with $t$ vertices. 
Let $W_n,\mu_n$
    and $\sigma_n$ be as above, and $Y_n:=(W_n-\mu_n)/\sigma_n$. Then
    \begin{align}\label{eq:treem}
  \mu_n&=\Theta(n),
\\\label{eq:treev}
\sigma^2_n&=O(n).              
    \end{align}
    If\/ $\sigma_n^2=\Theta(n)$, then there is a constant $C=C(m,t)$ such that 
    \begin{equation}\label{eq:treena}
        \dw(\law(Y_n),\mathcal{N}(0,1))\leq Cn^{-\frac{1}{2}
        }.
    \end{equation}
\end{theorem}

In the trivial cases $t=1$ and $t=2$ (a single vertex and a single edge,
respectively), $W_n$ is deterministic and $\sigma^2_n=0$; thus $Y_n$ is not
even defined.
We conjecture that these are the only cases where \eqref{eq:treena} does not
hold. 
\begin{conjecture}\label{Conjtree}
  If\/ $\gL$ is tree with at least 2 edges, then $\sigma^2_n=\Theta(n)$, 
and thus the normal approximation \eqref{eq:treena} holds.
\end{conjecture}
We show in \refS{SP2} that this holds at least for 
stars.
\begin{theorem}\label{th:na2star}
    Conjecture \ref{Conjtree} holds when $\Lambda$ is a star $\Sl$
with $\ell\ge2$ edges. 
\end{theorem}

The remaining subgraphs are the multicyclic ones.
%; for them we give only partial results, leaving extensions as an open problem.
We consider in particular the ones that have the following properties.
\begin{definition}\label{def:AB}
We say that a  graph is
\begin{itemize}

\item%\label{AAA} 
\emph{multicyclic} if it has at least two (not
necessarily edge- or vertex-disjoint) cycles; 

\item%\label{BBB} 
\emph{leaf-free} if it has no node of degree 1. 
\end{itemize}
\end{definition}

\begin{figure}
    \centering
    \includegraphics[scale=0.25]{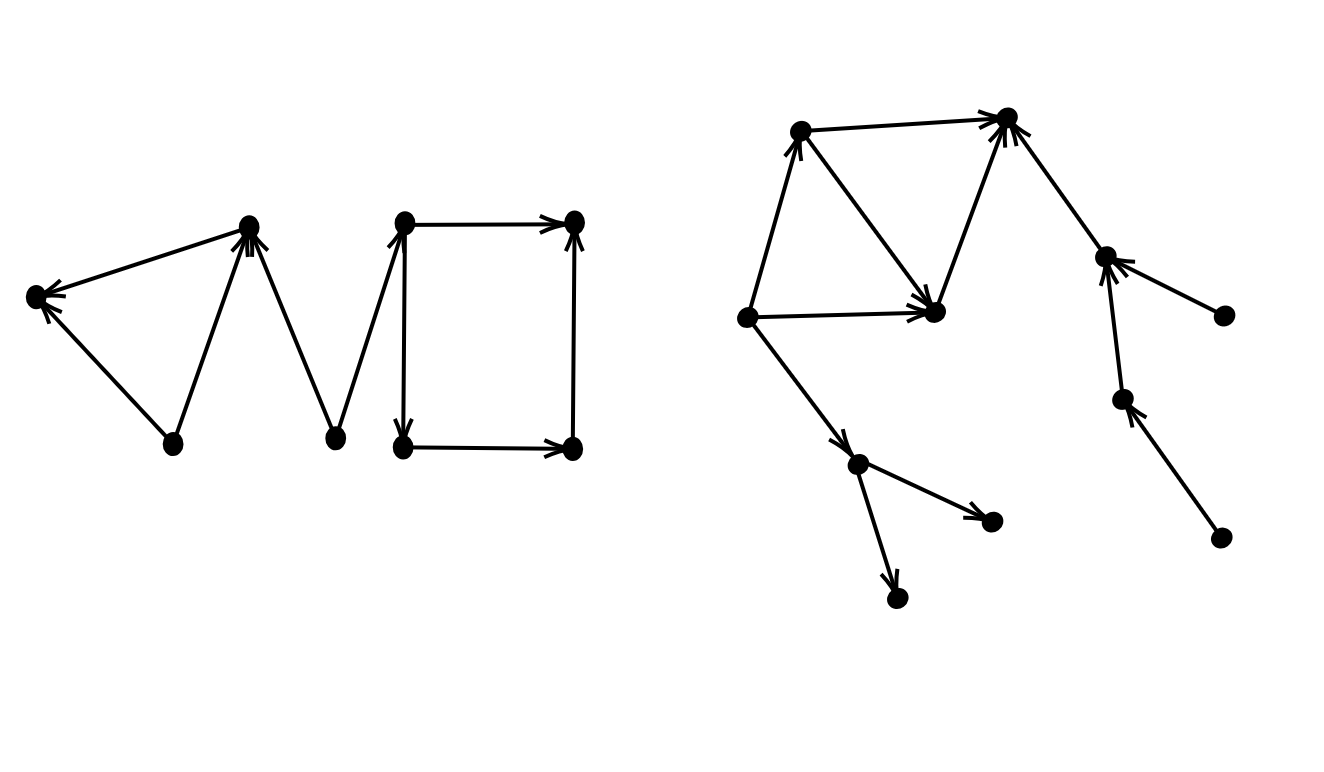}
 \caption{\small \emph{Left}: a graph that is both multicyclic and leaf-free. 
     \emph{Right}: a graph that is multicyclic but not leaf-free.
}
\label{fig:AB} 
\end{figure}

As examples of connected graphs that are both multicyclic and leaf-free, we
have two edge-disjoint cycles joined by an edge, or two cycles that share
precisely an edge; see Figure \ref{fig:AB} for another example. 

\begin{theorem}\label{th:bdexpsub}
    For any connected graph $H$ that is both multicyclic and leaf-free, and any $m\geq 2$,
    the expected number of copies 
    of $H$ in $G^m_n$ is bounded as $n\to\infty$.
\end{theorem}

\begin{remark}\label{RGoo}
  We may define the infinite random graph $G^m_\infty$ with vertex set
  \set{1,2,\dots} as the union of $G^m_n$ over all $n\ge1$.
Let $H$ be any fixed graph and let as above $W_n$ be the number of copies of
$H$ in $G^m_n$; also, let
$W_\infty\le\infty$ be the number of copies of
$H$ in $G^m_\infty$.
Then, as \ntoo, $W_n\upto W_\infty$.
In particular, $W_n\pto W_\infty$, and since both
\eqref{eq:ctm}--\eqref{eq:ctv} and \eqref{eq:treem}--\eqref{eq:treev} 
imply $W_n\pto+\infty$, we see that if $H$ is
unicyclic or a tree, then $W_\infty=\infty$ a.s.
On the other hand, 
\refT{th:bdexpsub} implies by monotone convergence that 
if $H$ is multicyclic and leaf-free, then 
$\E W_\infty<\infty$ and thus $W_\infty<\infty$ a.s.

Note further that since $W_n\asto W_\infty$ as \ntoo, in particular, $W_n$
converges in distribution to $W_\infty$ (without any normalization).
However, we do not expect that this limiting distribution has any nice form
such as Poisson, since $W_\infty$ is mainly determined by the
random wirings
of the first few edges in $G^m_n$.
\end{remark}

We also consider the expected counts of a graph that are multicyclic but not
leaf-free.
Note that every such  graph $H$ can be
constructed in the following way. Start with a graph $H'$ 
that is both multicyclic and leaf-free, and let $\cT_i$, $i=1,...,s$ be trees
with $t_i$
edges and a distinguished root. The graph $H$ is obtained by attaching
every $\cT_i$ to $H'$, with one of the vertices of $H'$ as the distinguished
root of $\cT_i$. As before, we may assume that any pair of $\cT_i$ do not share
the same vertex of $H'$. 
The graph $H'$ is known as the \emph{$2$-core} of $H$.
See Figure \ref{fig:AB} for an example. 

In the following theorem, we for convenience include the case $t=0$,
which is just \refT{th:bdexpsub}. (Then $s=0$ and $H'=H$.)
Recall the definition of $m_H$ in \refR{RmH}, and that $W_n=0$ if $m<m_H$.

\begin{theorem}\label{th:others}
    Fixing a connected multicyclic graph  $H$
and a positive integer $m\ge m_H$,
let $W_n$ be the
     number of copies of $H$ in $G^m_n$,
and let $t:=\sum_i t_i\ge0$ be the number of vertices in $H\setminus H'$,
where $H'$ as above is the $2$-core of $H$.
Then,
    \begin{align}
        \E W_n&=\Theta\bigpar{\logn{t}},\label{cq0}\\
        \var(W_n)&=O\bigpar{\logn{2t}}.\label{cq01}
    \end{align}
\end{theorem}
The precise order of $\var(W_n)$ is more difficult to
  establish for the same reason as for trees. For this family of graphs,
we do not expect the count to be 
approximately Poisson or normal (after renormalisation). 
Instead, we make the following conjecture, 
where we, as commented above, do not expect the distribution of the limit to
have a nice form. 
\begin{conjecture}\label{conj:A-B}
Let $H$ and $H'$ be as in \refT{th:others},
and let\/ $W_n$ and $W'_n$ be the numbers of copies of $H$ and $H'$ in
$G^m_n$, respectively.
Then there exists a constant $c>0$ such that
$(\log n)^{-t}W_n-c W'_n \asto 0$ as \ntoo, and thus
\begin{align}\label{cq1}
  \frac{W_n}{\log^t n} \asto c W'_\infty
.\end{align}
\end{conjecture}

We can prove a special case. It appears likely that the general case can be
shown by similar arguments, but the details seem complicated and we leave
this as an open problem.

\begin{theorem}\label{Ttail}
\refConj{conj:A-B} holds when $t=1$.
\end{theorem}

We note also that at least for some multicyclic graphs $H'$,
$W'_\infty=0$ with positive probability, and thus
$\P(W_\infty=0)\ge \P(W'_\infty=0)>0$; hence although $\E W_n\to\infty$ by
\eqref{cq0}, $\P(W_n=0)\ge c>0$, which in particular shows that we cannot
have Poisson or normal convergence.
We do not know whether this  holds for every simple multigraph $H'$,
and we give just one example.

\begin{example}
Let $H'$ be the complete graph on 4 vertices, $K_4$, and construct $H$ by
adding a single edge to $H'$; thus we take $s=1$ and $\cT_1$ as an edge in
the construction above. Let $W_n$ and $W'_n$ be as above.
By Theorem \ref{th:bdexpsub} (or a simple calculation), we see that
a.s.\ $W'_n\to W'_\infty<\infty$
as $n\to\infty$. 
Moreover, we claim that $\P(W'_\infty=0)>0$.

To see this, let $H''=K_4^-$ be $K_4$ minus one edge, and let $W_n''$ be the
number of copies of $K_4^-$ in $G^m_n$.
Fix a large $N$, and consider only $n>N$.
Let $\cE_N'$ be the event that there is a triangle in $G^m_N$.
Let $\cE_{N,n}''$ be the event that there is a copy of $K_4$ in $G^m_n\cup K^m_N$
with at most 2 vertices in $[N]$.
In other words $\cE''_{N,n}$ means that there is either a copy of $K_4$ in
$G^m_n$, or a copy of $K_4^-$ in $G^m_n$ with the two non-adjacent vertices
in $[N]$ and the two others in $[n]\setminus[N]$.
Note that $\cE_{N,n}''$ does not depend on the edges of $G^m_N$, and thus
$\cE_N'$ and $\cE_{N,n}''$ are independent.
If there is a copy of $K_4$ in $G^m_n$, then either $\cE_N'$ or $\cE_{N,n}''$
holds, and thus
\begin{align}
  \label{cq2}
\P(W_n'=0) \ge \P\xpar{\xnot\cE_N'}\P\xpar{\xnot \cE_{N,n}''}
.\end{align}
If $\cE_{N,n}''$ holds, then there is a copy of $H''=K_4^-$ in $G^m_n$ with at
least two vertices in $[n]\setminus[N]$, and thus
$W''_\infty \ge W''_n \ge W''_N+1$. Hence, by Markov's inequality,
\begin{align}\label{cq3}
  \P(\cE_{N,n}'') \le \E(W_\infty''-W_N'')
= \E W_\infty''-\E W_N''
.\end{align}
Since $\E W_N''\to \E W_\infty''$ as \Ntoo, we can choose $N$ such that
$\E W_\infty''-\E W_N''\le\frac12$; then
\eqref{cq3} implies that
$\P(\xnot\cE_{N,n}'')\ge \frac12$ for all $n>N$.
Moreover, $\P(\xnot\cE_N')>0$, since there is a positive probability that all
edges in $G^m_N$ lead to 1, and then there is no triangle.
Consequently, \eqref{cq2} yields $\P(W_n'=0)\ge c$, for some $c>0$ that does
not depend on $n>N$; thus also $\P(W'_\infty=0)\ge c>0$.

If there are no copies of $H'$, then there can be no copies of $H$, and thus
we conclude
$\P(W_n=0)\ge c>0$ for all $n$. Hence, although $\E W_n\to\infty$ by 
\eqref{cq0}, $W_n$ does not converge in probability
to $+\infty$, and in particular $W_n$ cannot be asymptotically normal.
\end{example}

\subsection{Discussion on possible future avenues}
An important direction for future work is to verify if Conjecture
\ref{Conjtree} holds, or in other words, to prove that the variance of the
count of any tree with at least 2 edges is of the precise order $n$.
Another direction for future work is to compute the leading
coefficients of the means and variances of the subgraph counts. More
precise expressions will possibly also enable us to verify if the approximation
results still hold if $m$ and (or) the number of vertices of the subgraph are
allowed to increase with $n$. 
It is also possible to prove a multivariate
analogue for the normal approximation results. This can be done using, for
instance, \cite[Theorem 12.1]{chen2010normal}, which also uses Stein's
method with the size-bias coupling.

\subsection{Article outline}
In the next section, we state the results from Stein's method that we apply
in the approximation proofs.
These results use a coupling that we construct 
in Section \ref{se:sbc}.
In Section \ref{se:prep}, we prove Theorem \ref{th:bdexpsub} and some
additional lemmas that will be useful in the approximation proofs
later. We prove the Poisson approximation results 
(\refTs{th:pocycle} and \ref{th:mpocycle})
in Section \ref{se:pocycles} and the normal approximation for 
unicyclic graphs (\refT{th:nact}) in
Section \ref{se:naballoons}, where we also prove \refT{th:others}. 
Section \ref{se:natrees} contains the proofs of
the normal approximation for  trees 
(\refTs{th:natree} and \ref{th:na2star}),
and in the last section we give a proof of \refT{Ttail}.

\section{Preliminary: Stein's method}\label{se:stein}
The error bounds in the Poisson and normal approximation results that we use
are obtained from general results on Stein's method
in terms of a coupling that we now describe. Let $(I_\al)_{\al\in
  \Gamma}$ be a collection of 0-1 valued random variables. For each
$\al\in\Gamma$, let the random variables $(J_{\beta\al})_{\beta\in\Gamma}$
be defined on the same space as $(I_\al)_{\al\in\Gamma}$, satisfying 
\begin{equation}\label{def:sbc}
    \law(J_{\beta\al}; \beta\in\Gamma\setminus\clc{\al})=\law(I_\beta;\beta\in \Gamma\setminus{\al}\mid I_\al=1).
\end{equation}
Note that this is a special case of the size-bias coupling appearing in the literature of Stein's method; see \textcite{chen2010normal}, \textcite{ross2011fundamentals}, and also \textcite{goldstein1996multivariate}. The theorem below is a direct consequence of  \textcite[Chapter 2, equation (1.2)]{barbour1992poisson}.

\begin{theorem}\label{th:pa}
Let $(I_\al)_{\al\in \Gamma}$ be as above, where $\E I_\al=\pi_\al$. Suppose that for each $\al\in\Gamma$, there is a coupling of $(J_{\beta\al})_{\beta\in\Gamma}$ and $(I_\al)_{\al\in\Gamma}$ such that \eqref{def:sbc} holds. Let $W:=\sum_{\al\in\Gamma}I_\al$, $\lambda:=\sum_{\al\in\Gamma}\pi_\al$. Then,
\begin{equation}\label{eq:poidtv}
    \dtv(\law(W),\mathrm{Po}(\lambda))\leq \min\clc{1,\lambda^{-1}}\bbclr{\sum_{\al\in\Gamma}\pi^2_\al +\sum_{\al\in\Gamma} \sum_{\beta\not=\al}\pi_\al \E | J_{\beta\al}-I_\beta|}.
\end{equation}
\end{theorem}

Now, assume that there is a partition $\Gamma=\bigcup^r_{j=1}\Gamma_j$. The next theorem quantifies how well the random vector $(\sum_{\al\in\Gamma_j}I_\al)_{1\leq j\leq r}$ can be approximated by a vector of independent Poisson variables. In our application, $\Gamma_j$ is to be taken as the set of potential cycles of length $\ell(j)$.
\begin{theorem}[{\cite[Theorem 10.K]{barbour1992poisson}}]\label{th:mpa}
Let $W_j=\sum_{\al\in\Gamma_j}I_\al$, $\E I_\al=\pi_\al$ and $\lambda_j=\E W_j$. Suppose that for each $\al\in\Gamma$, there is a coupling of $(J_{\beta\al})_{\beta\in\Gamma}$ and $(I_\al)_{\al\in\Gamma}$ such that \eqref{def:sbc} holds. Then,
\begin{equation}
    \dtv\bbclr{\law(\clc{W_j}^r_{j=1}), \prod^r_{j=1}\mathrm{Po}(\lambda_j)}\leq \frac{1+2\log^+(e\min_j \lambda_j)}{e\min_j\lambda_j}\sum_{\al\in\Gamma} \bbclr{\pi^2_\al + \sum_{\beta\not=\al}\pi_\al \E | J_{\beta\al}-I_\beta|}.
\end{equation}
\end{theorem}

We now state a normal approximation result for a collection of 0-1
variables, which follows from 
\cite[Theorem 3.20 and Corollary 3.24]{ross2011fundamentals}. 

\begin{theorem}\label{th:wassersteinbd}
Let $(I_\al)_{\al\in\Gamma}$ be a collection of $0$-$1$ variables and let
$(J_{\beta\al})_{\beta\in\Gamma}$ be defined on the same space as
$(I_\al)_{\al\in\Gamma}$, satisfying \eqref{def:sbc}. Define
$W:=\sum_{\al\in\Gamma}I_\al$, $\mu:=\E W$, $\sigma^2:=\var(W)$ and $W^s:= 
\sum_{\beta\in \Gamma\setminus \{K\}}J_{\beta K}+1$, where the index
$K\in\Gamma$ 
is chosen randomly with probabilities $\IP(K=\alpha)=\E I_\alpha/\mu$. 
If $Z=(W-\mu)/\sigma$, then
\begin{equation}\label{eq:wassersteinbd}
    \dw\bclr{\law(Z),\mathcal{N}(0,1)}\leq \frac{\mu}{\sigma^2}\sqrt{\frac{2}{\pi}}\sqrt{\var\bclr{\E[W^s-W\mid W]}} + \frac{\mu}{\sigma^3}\E[(W^s-W)^2],
\end{equation}
where $\mathcal{N}(0,1)$ is the standard normal distribution.
\end{theorem}

\section{Construction of the size-bias couplings}\label{se:sbc}
As the edges of a uniform attachment graph are independent, we can construct
the coupling appearing in Theorems \ref{th:pa}--\ref{th:wassersteinbd} 
%and \ref{th:mpa} 
as
follows. Fix the subgraph $H$, let
$G:=G^m_n$ be the uniform attachment graph and let $\Gamma$ be the set of
potential copies of $H$. For every $\al\in \Gamma$, we couple two graphs
$G$ and $G^\al$ by matching their attachment steps, except for the edges of~$\al$. In $G^\al$, the edges of $\al$ are wired in a deterministic fashion
to obtain the copy $\al$ of $H$; whereas in $G$, they are generated independently
from the construction of $G^\al$. 
For $\gb\in\Gamma$,
let as above $\tone_\beta$ 
be the indicator that the subgraph $\gb$ is present in $G$, and 
let $\tone^\ga_\beta$ 
be the indicator that $\gb$ is present in $G^\ga$.
It follows that for every chosen $\alpha$, 
\begin{equation}\label{1gagb}
\law(\tone^\al_\beta;\beta\in\Gamma\setminus\clc{\al})=\law(\tone_\beta; \beta\in\Gamma\setminus\clc{\al}\mid\tone_\al=1).
\end{equation}
Let $\Gamma^-_\al\subset\gG$ be the set of copies $\gb$ of $H$ such that at
least 
one edge of $\al$ has a different endpoint in $\gb$, and let
$\Gamma^+_\al\subset\Gamma\setminus \Gamma^-_\al$ be the set of copies of
$H$ that share at least one edge with $\al$, excluding $\al$ itself. Observe
that under the coupling, 
$\tone^\al_{\beta}\leq \tone_\beta$ for $\beta\in\Gamma^-_\al$
and $\tone^\al_{\beta}\geq \tone_\beta$ for $\beta\in\Gamma^+_\al$,
while $\tone^\al_{\beta}=\tone_\beta$ for every
$\beta\in\Gamma\setminus(\Gamma^+_\al\cup \Gamma^-_\al\cup \clc{\al})$. The
error bound (\ref{eq:poidtv}) in Theorem \ref{th:pa} therefore simplifies to 
\begin{equation}\label{eq:simplifiedsteinbd}
    \min\clc{1,\lambda^{-1}}\lrpar{\sum_{\al\in\Gamma}\pi_\al^2 + \sum_{\al\in\Gamma}\sum_{\beta\in \Gamma^+_\al} \pi_\al\E[\tone^\al_{\beta}-\tone_\beta] + \sum_{\al\in\Gamma}\sum_{\beta\in \Gamma^-_\al} \pi_\al\E[\tone_\beta-\tone^\al_{\beta}]};
\end{equation}
noting that $\sum_{\beta\in\Gamma^+_\al}\E[\tone^\al_\beta-\tone_\beta]$ and
$\sum_{\beta\in \Gamma^-_\al} \Pa\E[\tone_\beta-\tone^\al_{\beta}]$ are
respectively the expected gain and loss in the copies of $H$ after forcing
$\al$ to be present in the graph $G$. 
We may simplify this further by omitting the negative terms in
\eqref{eq:simplifiedsteinbd} and noting that, by \eqref{1gagb},
\begin{align}
  \pi_\ga\E\lrsqpar{\tone^\ga_\gb}
=  \P(\tone_\ga=1)\E\lrsqpar{\tone_\gb\mid \tone_\ga=1}
%=\P\lrsqpar{\tone_\gb \tone_\ga=1}
=\E\lrsqpar{\tone_\gb \tone_\ga}
,\end{align}
which yields the simpler (but somewhat larger) error bound
\begin{equation}\label{eq:simplifiedsteinbd+}
    \min\clc{1,\lambda^{-1}}\lrpar{\sum_{\al\in\Gamma}\pi_\al^2 
      + \sum_{\al\in\Gamma}\sum_{\beta\in \Gamma^+_\al}\E\lrsqpar{\tone_\ga \tone_\gb}
      + \sum_{\al\in\Gamma}\sum_{\beta\in \Gamma^-_\al} \pi_\al\pi_\beta}.
\end{equation}

As for Theorem~\ref{th:mpa}, the simplified error bound is the same as
\eqref{eq:simplifiedsteinbd} or \eqref{eq:simplifiedsteinbd+}, 
but with the factor $\min\{1,\lambda^{-1}\}$ replaced with 
\begin{equation}
    \frac{1+2\log^+(e \min_j \lambda_j )}{e\min_j \lambda_j}
\end{equation}

When applying Theorem \ref{th:wassersteinbd}, we first sample a copy $K\in
\Gamma$ of subgraph $H$ with probabilities $\P(K=\ga)$ proportional to
$\E\tone_\ga$, 
%proportional to its mean,
and construct the graphs $G^K$ and $G$ as above. The subgraph count $W$ and its
size-bias version $W^s$ can then be found in $G$ and $G^K$, respectively.

%%XXXXXXXXXXXXXXXXXXXXXXXXXXXXXXXXXXXXXXXXXXXXXXXXXXXXXXXXXXX

%\input{preplemmas}
\section{Proof of Theorem \ref{th:bdexpsub} and some 
useful lemmas}\label{se:prep}

Here we prove some lemmas that are useful for proving the main results
later. We also prove Theorem \ref{th:bdexpsub} here, as some special cases
of the result will be applied in the other proofs. 
To study the expeted number of copies of a graph  $H$, we use the following
definition. 
Let $h$ be the number of vertices in $H$.

\begin{definition}[Vertex marks and mark sequence]\label{def:marks}
Let $\ga$ be a potential copy of $H$ and suppose that 
its vertices are $k_1<\dots <k_h$.
We say that the
\emph{mark} of a vertex $k_i\in \ga$ is the out-degree of $k_i$ in 
$\ga$ regarded as a directed graph (as always, with edges directed towards
the smaller endpoint).
The \emph{mark sequence} is $(b_i)_{1\leq i\leq h}$, with $b_i$ being the
mark of vertex~$k_i$. 
\end{definition}

In other words, for $\ga$ to actually be a copy of $H$ in $G^m_n$, 
there are for each $i\in[h]$ exactly $b_i$ edges in $G^m_n$ from $k_i$ 
that have to have the endpoints determined by $\ga$.
Note that the mark sequence does not entirely encode the configuration of
the copy of $H$ in $G^m_n$, 
but, together with the sequence of vertices $k_i$,
the mark sequence gives the probability that the copy $\ga$ is
present in $G^m_n$. In fact, since the edges in $G^m_n$ choose their
endpoints independently,
\begin{align}\label{pga}
\Pa:=
  \P(\tone_\ga=1) = \prod_{i=1}^h \frac{1}{(k_i-1)^{b_i}}
= \prod_{i=2}^h \frac{1}{(k_i-1)^{b_i}}
,\end{align}
where the final equality holds because always $b_1=0$, since edges are
directed towards lower
numbers, and thus  $\ga$ has no out-edge from $k_1$. Similarly,
$b_h$ equals the degree of $k_h$ in $\ga$, since $k_h$ has no in-edge.

Note that for a given unlabelled graph $H$, there is a finite number 
(at most $h!$) of non-isomorphic labelled versions of $H$.
Hence, in order to obtain estimates of the expected number of
copies, it suffices to consider each possible labelled version separately,  
and then consider only potential copies with vertex sequences $(k_i)_1^h$ where
$k_1<\dots<k_h$  and $k_i$ corresponds to vertex $i\in H$.
The mark sequence of the potential copy 
depends only on the labelled version of $H$.
Different labelled versions of $H$ may yield different mark sequences
$(b_i)_1^h$,  
but they all
have the same length $h=v(H)$,  the number of vertices in $H$, and the same sum
\begin{align}\label{marksum}
  \sum_{i=1}^h b_i=e(H),
\end{align}
the number of edges in $H$.

For the proof of Theorem \ref{th:bdexpsub}, we also use the following
definition. 

\begin{definition}[$F$-number] \label{def:Fnumber}
 We define the \emph{$F$-number}, $F(a_i)$ of a finite sequence $(a_i)$ of natural numbers as $\sum_{a_i>1}a_i+\sum_{a_i=0}(-2)$, that is, we ignore all $a_i$ with $a_i=1$, count the $a_i$ with $a_i=0$ with weight $-2$ and count all other $a_i$ with $a_i$ as their weight.
\end{definition} 

\begin{proof}[Proof of Theorem \ref{th:bdexpsub}]
Consider a connected graph $H$ that is multicyclic and leaf-free,
see
Definition \ref{def:AB}. 
Let $(b_i)_{1\leq i\leq h}$ be the mark sequence of a potential copy of
$H$.
By \eqref{pga}, the expected number of copies of $H$ with this mark sequence can
be bounded by
a constant times
\begin{equation}\label{eq:expABH}
S(b_1\xdots b_h):=
    \sum_{1\le k_1<\dots<k_h<\infty}\prod^h_{i=2}\frac{1}{(k_i-1)^{b_i}}.
\end{equation}
Hence it suffices to show that this infinite sum is finite
for every mark sequence $(b_i)_1^h$. 

We will show this by modifying the mark sequence $(b_i)$ in ways that can
only increase the sum $S(b_1\xdots b_h)$ in \eqref{eq:expABH}, until we
reach a sequence where we can show that $S(b_1\xdots b_h)$ is finite.
Note that we do not claim that the modified sequences actually are mark
sequences for some copies of $H$; we consider the value 
$S(b_1\xdots b_h)$ as defined by \eqref{eq:expABH}, without worrying about a
probabilistic interpretation in general.

Since $H$ is leaf-free, each vertex in $H$ has degree at least 2, and thus
$b_h>1$. 
We change the sequence $(b_i)_1^h$ as follows: 
In the first round,
if $b_h>2$ 
and $b_i=0$ for some $i\neq1$,
then decrease $b_h$ by 1 and
increase the last such $b_i$ (i.e., the one with maximal $i$) by 1.
For instance, $0102023 \rightarrow  0102122$. We
repeat this process until $b_h=2$ or $b_i>0$ for all $i\neq 1$. 
In the second round, if
$b_h=2$ and $b_{h-1}>1$, then 
we repeat the same procedure with $b_{h-1}$ until $b_{h-1}=1$ or we
have exhausted all $b_i$, $i\neq1$, such that $b_i=0$. 
Continue in the same way: in the $p$th round ($3\le p<h$), 
if $b_h=2$, $b_{h-1}=\dots=b_{h-p+2}=1$, and $b_{h-p+1}>1$,  
then repeat the
procedure for $b_{h-p+1}$ as we did for $b_{h-1}$. Stop when
such moves are no longer possible. 
Denote the final sequence by $(b'_i)_1^h$.
Note that these moves never decrease
$\prod^h_{i=1}(k_i-1)^{-b_i}$ since we always shift mass in the mark sequence
to the left; 
thus $S(b_1\xdots b_h)\le S(b'_1\xdots b'_h)$. 

All the possible final sequences $(b'_i)_1^h$ are analysed below:
\begin{enumerate}
\item 
We have at  least one intermediate zero left, say $b'_\jq=0$ with $1<\jq<h$
and $\jq$ maximal among such indices. 
Then $b'_{\jq+1}=\dots=b'_{h-1}=1$ and $b'_h=2$, since otherwise we would have
modified the sequence further. In other words, 
we have a sequence of the form $0XXX\xdots XXX011\xdots 112$. 
This is not possible, due to the $F$-number of the final part of the sequence
$b'_\jq\xdots b'_h=0111\xdots 1112$. Note that the
moves that we made above cannot have affected $b_i$ for $i<\jq$, so they
were entirely in this subsequence, and it is easy to see that each move
either preserved or increased the $F$-number of this
subsequence. Thus $0=F(0111\xdots 1112)\geq F(b_{\jq}\xdots b_h)$.
with $(b_i)_{1\leq i\leq h}$ being the original mark sequence. A way of
interpreting the $F$-number is as an upper bound of the number
of edges directed from the
last $h-\jq+1$ nodes to the first $\jq-1$ nodes,
which equals the sum of the differences between outdegree and indegree for
the last $h-\jq+1$ nodes. 
In any copy of a connected graph
$H$ that is multicyclic and leaf-free, a vertex with mark $0$ must receive at least
$2$ in-edges and thus contributes at most $a_i=-2$ to the count. A vertex
$i$ with mark $1$ sends one edge but receives at least $1$, and so
contributes at most $a_i=0$ to the count. Similarly, a vertex $i$ with mark
$k>1$ may receive some
in-edges and so contributes at most $a_i=k$ to the count. Since $H$ is a
connected graph, we must have some edge connecting the last $h-\jq+1$ vertices
and the first $\jq-1$ vertices in $H$, and so $F( b_{\jq}\xdots b_h)>0$,
contradicting the above.

    \item We have removed all the zeroes (except the first) and ended up
      with $0XX\xdots XX b'_h$ where $b'_h>2$ and each $X$ is at least
      $1$. (In this case, we could not reduce the last number down to 2
      before running out of zeroes.) The sequence can be compared to the
      sequence $0111\xdots 1113$ which gives a larger value of
      (\ref{eq:expABH}). For the latter sequence we have, if $h\ge3$,
by first summing over $k_h>k_{h-1}$,
\begin{align}\label{aq2}
S(011\xdots113)&=
\sum_{1\le
  k_1<\dots<k_h}\prod^{h-1}_{i=2}\frac{1}{k_i-1}\cdot\frac{1}{(k_h-1)^3}
\notag\\&
\le
\sum_{1\le k_1<\dots<k_{h-1}}\prod^{h-1}_{i=2}\frac{1}{k_i-1}\cdot\frac{1}{(k_{h-1}-1)^3}
,\end{align}
which is the same sum with $h$ replaced by $h-1$; thus induction yields
\begin{align}\label{aq3}
S(011\xdots113)\le 
S(03)=
\sum_{1\le k_1<k_2}\frac{1}{(k_2-1)^3}
= \sum_{k_2\ge2}\frac{1}{(k_2-1)^2}<\infty
.\end{align}

    \item We have $0XXX\xdots XXX2$ where the $X$'s are at least 1. In this
      case, note that not all the $X$'s can be 1 for topological reasons, as
      two connected cycles require more edges than nodes. Thus at least one
      $X$ is at least $2$. 
If $b'_2=1$, we can exchange $b'_2$ with $b'_j$ for some $j<h$ such that
$b'_j>1$, again without decreasing $S$.  
We can then compare the sequence to $02111\xdots 1112$.
By arguing as in \eqref{aq2}--\eqref{aq3} we find
\begin{align}\label{aq4}
  S(0211\xdots112)&
\le S(022)
=\sum_{1\le k_1<k_2<k_3}\frac{1}{(k_2-1)^2(k_3-1)^2}
\notag\\&
\le \sum_{1\le k_1<k_2}\frac{1}{(k_2-1)^3}=S(03)<\infty.
\end{align}
\end{enumerate}
This completes the proof of the theorem.
\end{proof}

For the expected number of $\ell$-cycles, 
which is not covered by \refT{th:bdexpsub}, 
we use a simpler version of the argument above to
prove the following lemma,
%the main ingredient is the following lemma, 
which essentially says that it is enough to just consider
one particular configuration of an $\ell$-cycle. We first observe that in 
the mark sequence of
a potential copy of a cycle, 
all marks are 0, 1, or 2; furthermore,
there must be equal numbers of marks  0 and 2 by \eqref{marksum}.  
An $\ell$-cycle has $\ell$ marks. 
As noted above, we must have $b_1=0$ and $b_\ell=2$,
since edges are directed towards lower numbers.
Note that an $\ell$-cycle has at most $\floor{\ell/2}$ vertices of mark 0.

\begin{lemma}\label{le:cfacts}
Among all potential copies of an $\ell$-cycle on given vertices
$k_1<\dots<k_\ell$, the configurations with the mark sequence $011\xdots112$ have the largest probability
to occur in~$G^m_n$.
\end{lemma}

\begin{proof} 
%  For each vertex $i\in[\ell]$, consider its in- and out-edges that are part
%  of the $\ell$-cycle. 
In a cycle, if a vertex sends $j$ edges, it receives
  $2-j$ edges. 
%It follows directly that the last vertex $\ell$ has mark 2
%  and vertex 1 has mark 0. 
Consider the mark sequence $(b_i)_1^\ell$ of a
  potential copy of an $\ell$-cycle. 
Counting from the right, and considering the last $k$
  vertices $\ell-k+1,\dots, \ell$, the number of marks 2 must
  always be larger than the number of marks 0, except when
  $k=\ell$, in which case they are equal;
equivalently, 
\begin{align}\label{aq1}
\sum_{\ell-k+1}^\ell (b_i-1)\ge 0, 
\text{with strict inequality when $1\le k<\ell$}.  
\end{align}
This holds since the cycle
  is connected, and so the total number of out-edges of the last $k$
  vertices must be larger than their total number of in-edges when $k<\ell$
  ($k=\ell$ gives equality). The probability of the cycle appearing in $G^m_n$
  is given by \eqref{pga}.
Consider all sequences $(b_i)_1^\ell$ with $b_i\in \set{0,1,2}$ such that
\eqref{aq1} holds. We must have $b_1=0$ and $b_\ell=2$. If there is another
mark $b_j=2$ with  $j\neq\ell$
in the sequence, let $j$ be the largest such index.
Then there must also be an index $i\neq1$ with $b_i=0$; let $i$ be the
largest such index. We must have $i<j$, since otherwise \eqref{aq1} would be
violated for $k=\ell-i+1$.
It is then easy to see that if we change both $b_i$ and
$b_j$ to 1, then \eqref{aq1} still holds, and the value of \eqref{pga} is
increased.
Consequently, the mark sequence with the largest probability is of the form $0111\xdots 1112$.
\end{proof}

The lemma below concerns  trees, and is used when we prove the normal
approximation results.

\begin{lemma}\label{le:gentree}
Fix the positive integer $t$ and let $\cT$ be a rooted tree with $t$
edges and thus $t+1$ vertices.
The following results concerning $\cT$ hold.
\begin{enumerate}[label=\textup{(\roman*)}]
    \item\label{le:gentree1}  
Among all potential copies of $\cT$ on a given vertex set, 
the configurations with the mark sequence $011\xdots11$ have the largest
probability.
    \item\label{le:gentree2} 
Let $N_x$ be the number of copies of
$\cT$ with vertex $x$ as its distinguished root. Then 
$\E N_x=O\xpar{\logn{t}}$,
uniformly for $1\le x\le n<\infty$.
\end{enumerate}
\end{lemma}

\begin{proof}
\ref{le:gentree1}:
Let $0b_2\xdots b_{t+1}$ be the mark sequence of a potential copy of $\cT$.
Note first that if $2\le k\le t+1$, then we have, summing over the $t+2-k$ last
vertices,
\begin{align}\label{dq1}
  \sum_{i=k}^{t+1}b_i\ge t+2-k.
\end{align}
In fact, this sum is the number of edges with upper endpoint in $[k,t+1]$,
which equals the number of edges with at least one endpoint in $[k,t+1]$,
and in any connected graph, a proper subset 
of $\ell$ vertices is always adjacent to at least $\ell$ edges. (To see
this, collapse all other vertices into one, 
and pick a spanning tree with $\ell$ edges  in the resulting graph.)

We argue similarly to the proof of Theorem \ref{th:bdexpsub},
now modifing the mark sequence as follows.
As long as some $b_i\ge2$, we reduce the rightmost such $b_i$ by 1 and
increase the rightmost mark 0 to 1. 
We can never reach a sequence ending with a proper subsequence
$011\xdots11$ with some $b_k=0$ ($k\ge2$) and $b_{k+1}=\dots=b_{t+1}=1$,
because the first time this happens, all previous moves must have been
inside $[k,\dots,t+1]$, so the sum in \eqref{dq1} has not changed
and \eqref{dq1} still holds, a contradiction.
Consequently, we do not stop until there is no mark 0 except $b_1$, but
since the sum of all marks is $t$ by \eqref{marksum} (and this sum is not
changed), the final sequence is $011\xdots11$.
As the procedure never decreases the probability \eqref{pga},
the mark sequence $011\xdots11$ indeed yields
the largest probability. 

\ref{le:gentree2}:
We use \ref{le:gentree1} to obtain
\begin{align}  
\E N_x&
\leq C\lrpar{\sum_{x<k_1<\dots<k_{t}\le n}\prod^{t}_{i=1}\frac{1}{k_i}
+ \sum^{t}_{r=1}\sum_{k_1<\dots<k_r<x<k_{r+1}<\dots<k_{t}\le n}
 \frac{1}{x}\prod^{t}_{i=2}\frac{1}{k_i}}
\notag\\&
\leq C\lrpar{\sum_{x<k_1<\dots<k_{t}\le n}\prod^{t}_{i=1}\frac{1}{k_i}
+ \sum^{t}_{r=1}\sum_{k_1<\dots<k_r<x<k_{r+1}<\dots<k_{t}\le n}
\prod^{t}_{i=1}\frac{1}{k_i}}
\notag\\&
\le C \prod^{t}_{i=1}\sum_{k_i=1}^n\frac{1}{k_i}
\le C\bigpar{\log n+1}^t
,\end{align}
which completes the proof. 
\end{proof}

%%XXXXXXXXXXXXXXXXXXXXXXXXXXXXXXXXXXXXXXXXXXXXXXXXXXXXXXXXXXX
%%\input{cycles}

\section{Proof of the Poisson approximations for cycles}\label{se:pocycles}
In this section we prove the Poisson approximation results in Theorems
\ref{th:pocycle} and \ref{th:mpocycle} using 
Theorems \ref{th:bdexpsub}, \ref{th:pa} and \ref{th:mpa}. 
We fix in this section the  integers $m\ge1$ and $\ell\ge2$, 
and denote by $W_n$ the count of
$\ell$-cycles in $G^m_n$.

\begin{figure}
    \centering
    \includegraphics[scale=0.4]{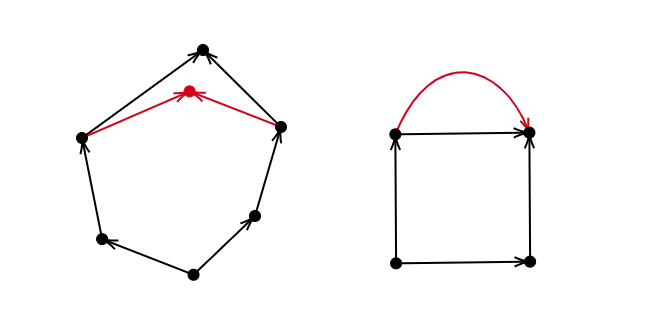}
    \caption{{\small Two examples of positively correlated cycles. In each example, the red edges form the additional segment that gives rise to the other cycle.}}
    \label{fig:pos-cycles}
\end{figure}

\begin{proof}[Proof of Theorems \ref{th:pocycle} and \ref{th:mpocycle}]
We start by proving the result on the expected number of cycles $\mu_n$ in
(\ref{eq:cyclem}). 
%Fix the parameters $\ell$ and $m$. 
By Lemma
\ref{le:cfacts} and \eqref{pga}, $\mu_n$ can be bounded above by 
\begin{align}
\psi_{\ell,m}\sum_{k_1<\dots<k_\ell\le n}\frac{1}{(k_\ell-1)^2}\prod^{\ell-1}_{i=2}\frac{1}{k_i-1},
\end{align}
where $\psi_{\ell,m}$ is the total number of ways of forming a potential
$\ell$-cycle on $\ell$ given vertices. Summing over
$k_1,\dots,k_{\ell-1}$ (in this order), the above is at most
$\psi_{\ell,m}\sum^n_{k_\ell=\ell}\frac{1}{k_\ell-1}$, which can be bounded
by $\psi_{\ell,m} (\log n+1)$. The calculation for the lower bound is similar,
with the factor $\psi_{\ell,m}$ replaced by 1; we consider just one
configuration that has mark sequence $011\xdots11$.

Next, we 
construct the size-bias coupling of the cycle count described in Section
\ref{se:sbc}.
To prove the approximation result in Theorem \ref{th:pocycle}, it then
suffices to show that the sums in (\ref{eq:simplifiedsteinbd+}) are
bounded as $n\to\infty$; the error bound then follows from
\eqref{eq:simplifiedsteinbd+} and
(\ref{eq:cyclem}). 

Let $\Gamma$ be the set of potential cycles and let as in \eqref{pga} $\Pa$
be the probability of cycle~$\al\in\Gamma$ being a subgraph of $G^m_n$. 
By \eqref{pga} and \refL{le:cfacts}, using the notation \eqref{eq:expABH},
\begin{align}\label{fq0}
\sum_{\al\in
  \Gamma}\Pa^2\leq
C\sum_{k_1<\dots<k_\ell\le n}k^{-4}_\ell\prod^{\ell-1}_{i=2}k^{-2}_i
\le C\cdot S(022\xdots224)
\le C \cdot S(011\xdots113),
\end{align}
which is finite by \eqref{aq3}.

In the sum
\begin{align}\label{fq1}
\sum_{\al\in\Gamma}\sum_{\beta\in\Gamma^+_\al}
\E\lrsqpar{\tone_\ga \tone_\gb}
, 
\end{align}
each term $\E\sqpar{\tone_\ga \tone_\gb}=\P\sqpar{\tone_\ga \tone_\gb=1}$
is the probability that $G^m_n$ contains a specific copy of a graph $H$
obtained by merging the cycles $\ga$ and $\gb$.
See Figure \ref{fig:pos-cycles} for an illustration.
There is only a finite number of such graphs $H$
(up to isomorphism); each copy of a graph $H$
arises for a bounded number of pairs $(\ga,\gb)$; and
each such graph $H$ is 
connected, 
multicyclic and leaf-free (see Definition \ref{def:AB}).
Consequently, the sum \eqref{fq1} is bounded by a constant times the sum of
the expected number of copies of such graphs $H$, and 
Theorem \ref{th:bdexpsub} shows that this sum
is bounded as $n\to\infty$.

Finally, note that if $\gb\in\gG^-_\ga$, then the cycles $\ga$ and $\gb$
have at least one common vertex. 
%Fix $\al\in\Gamma$ and  $\beta\in \gG^-_\ga$. 
Consider a uniform attachment graph $G^{2m}_n$ where
  the first $m$ out-edges of any vertex $2\leq i\leq n$ are coloured red,
  and the remaining $m$ out-edges of these vertices are coloured blue. 
This graph can be thought of as overlaying two independent copies
of $G^{m}_n$. Let $\gaxb$ be the subgraph of $K^{2m}_n$ obtained by
regarding $\al$ and $\beta$ living in the red and blue parts respectively,
and taking the union of them. 
(Thus, the $i$-th edge of vertex~$j$, $j^{(i)}$, that is part of $\beta$ is
now written as $j^{(i+m)}$.) 
With some slight abuse of notation, let
$\tone_{\gaxb }$ be the indicator of the subgraph $\gaxb$ in $G^{2m}_n$. By
the edge independence, we have $\pi_\al\pi_{\beta}= \E\tone_{\gaxb}$,
which is the probability that $G^{2m}_n$ contains a given copy of
a graph $H$ that, just as in the argument for $\gG ^+_\ga$,
is obtained by merging the cycles $\ga$ and $\gb$.
There is only a finite number of such graphs $H$, so it follows from
Theorem \ref{th:bdexpsub} that
the expected number of copies of any of them in
$G^{2m}_n$ is bounded
as $n\to\infty$, and thus
\begin{align}\label{fq2}
\sum_{\al\in\Gamma}\sum_{\beta\in \Gamma^-_\al}\Pa\Pb 
=\E \sum_{\al\in\Gamma}\sum_{\beta\in \Gamma^-_\al}\tone_{\gaxb}
=O(1).
\end{align}
This completes the proof of \refT{th:pocycle}.

The proof of Theorem \ref{th:mpocycle} is similar, using \refT{th:mpa}; 
we now have
$\gG=\bigcup_{j=1}^r\gG_j$, where $\gG_j$ is the set of potential
$\ell(j)$-cycles. We estimate \eqref{fq0}--\eqref{fq2} as above,
when necessary
replacing $\ell$ by $\ell(j)$ and summing over $j=1,\dots,r$. 
\end{proof}

%%XXXXXXXXXXXXXXXXXXXXXXXXXXXXXXXXXXXXXXXXXXXXXXXXXXXXXXXXXXX
%%\input{balloons}

\section{Proof of the normal approximation for unicyclic
  graphs}\label{se:naballoons}

\subsection{Proofs of \eqref{eq:ctm} and \eqref{cq0} }\label{Sctm}
We first prove \eqref{eq:ctm} in \refT{th:nact} and \eqref{cq0} in \refT{th:others}, which
have similar proofs.

\begin{proof}[Proof of \eqref{eq:ctm}]
%We first prove the statement on $\mu_n$. 
Given a potential $\ell$-cycle 
$\al$, let
$\Xi(\al)$ be the set of 
$s$-tuples $(\gb_1,\dots,\gb_s)$ of
potential copies of $\cT_i$, $i=1,\dots,s$, that 
can be added to $\ga$ to form a potential copy of $\gL$.
Denote by $\Gamma_{\mathrm{cycle}}$ the set of potential $\ell$-cycles. By
independence of the edges, (\ref{eq:cyclem}) and Lemma~\ref{le:gentree}(ii),
we have
\begin{align}\label{eq:exballoons}  
\mu_n
=\sum_{\al\in\Gamma_{\mathrm{cycle}}} \sum_{\beta\in \Xi(\al)} \E\tone_\al  \E \tone_\beta 
    \leq O\bigpar{\logn{t}}\sum_{\al\in\Gamma_{\mathrm{cycle}}}\E\tone_\al  
=O\bigpar{\logn{t+1}},
\end{align}
noting that in the inequality we have used the assumption that $\ell$ and
$\cT_i$, $i=1,\dots,s$, are fixed.

For the lower bound on $\mu_n$, we choose a suitable configuration on a set
of vertices $k_1<\dots<k_{\ell+t}$. We use $k_1<\dots<k_\ell$ to form a
cycle whose configuration has the mark sequence $011\xdots112$, 
and $k_{\ell+1}<\dots<k_{\ell+t}$ to construct the trees
$(\cT_i)_{1\leq i\leq s}$, taking their vertices in some order going from
the roots outwards, so that each non-root vertex in $\cT_i$ has mark 1; this gives a potential copy of $\gL$ with all marks 1
except $b_1=0$ and $b_\ell=2$. 
Hence $\mu_n$ is at least, for some constant $c>0$,
by first summing
iteratively over $k_1<\dots<k_{\ell-1}$,
\begin{align*}\label{fq3}
\sum_{k_1<\dots<k_{\ell+t}\le n}\frac{1}{(k_2-1)\cdots(k_{\ell-1}-1)(k_\ell-1)^2}\prod^{\ell+t}_{i=\ell+1} \frac{1}{k_i-1} 
&\ge c\sum_{\ell\le k_\ell<\dots<k_{\ell+t}\le n} \prod^{\ell+t}_{i=\ell} \frac{1}{k_i-1}\\ &=\Theta\bigpar{\logn{t+1}},\numberthis
\end{align*} 
as required.
\end{proof}

\begin{proof}[Proof of \eqref{cq0} in \refT{th:others}]
The upper bound
follows as in (\ref{eq:exballoons}), now letting $\ga$ be a potential
copy of $H'$ and summing over all such $\ga$, using
Theorem \ref{th:bdexpsub} and, as above, Lemma~\ref{le:gentree}(ii).
For the lower bound, it suffices to consider one fixed potential copy of
$H'$.
\end{proof}

\subsection{Proofs of variance estimates}

We next prove the variance estimates \eqref{eq:ctv}  in Theorem \ref{th:nact} 
and \eqref{cq01} 
in \refT{th:others}. 

\begin{figure}
    \centering
    \includegraphics[scale=0.3]{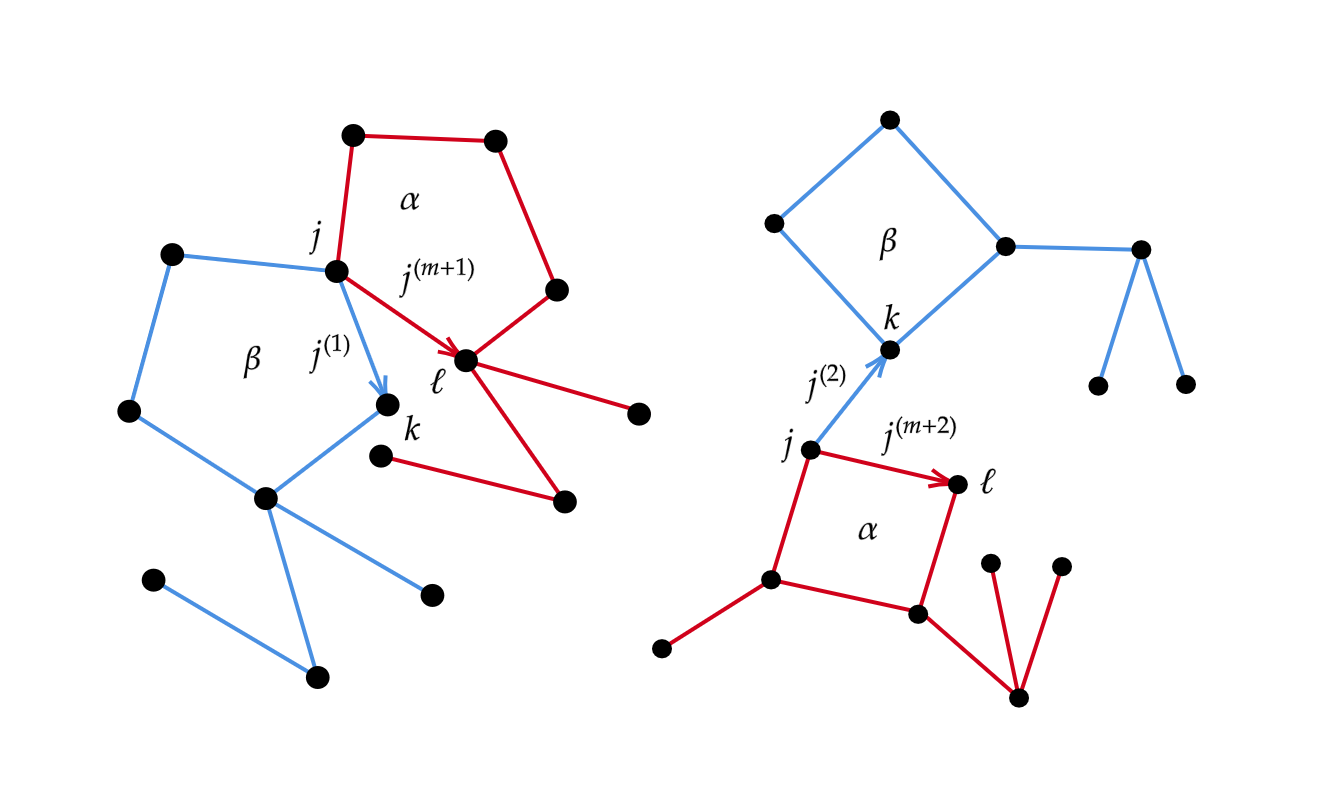}
    \caption{\small Two examples of two negative correlated copies $\al$ and $\beta$ that are embedded in the graph $G^{2m}_n$, where vertex $j$ is the common vertex and $j^{(1)}$ (left) and $j^{(2)}$ (right) are the edges that have different recipients in $\al$ and $\beta$.}
    \label{fig:cct}
\end{figure}

\begin{proof}[Proof of \eqref{eq:ctv}]
Recall that a connected unicyclic
graph $\Lambda$ can be constructed as follows. Let $\mathcal{C}_\ell$ be an
$\ell$-cycle and let $\cT_i$, $i=1,\dots,s$, be a tree with $t_i$ edges and a
distinguished root, so that $\cT_i$ has $t_i+1$ vertices. We construct
$\Lambda=\Lambda_{\ell,t_1,\dots,t_s}$ by attaching $\cT_i$ to
$\mathcal{C}_\ell$, using a vertex of $\mathcal{C}_\ell$ as the
distinguished root of $\cT_i$, and we may assume that $s\leq \ell$ and 
that at most one $\cT_i$ attaches to each vertex in $\mathcal{C}_\ell$. 
Fixing
$\Lambda$ throughout this section, we write
$\Gamma=\Gamma^{(n,m,\ell,t)}_{\Lambda}$ as the set of all potential copies
of $\Lambda$ in the uniform attachment graph $G:=G^m_{n}$. 

Given $\ga\in\gG$, denote by $\Gamma^-_\al$ the set of potential copies of
$\Lambda$ that
cannot coexist with $\al$ in the same instance of $G^m_n$ 
(because at least one edge is incompatible with the edges of $\al$), 
and denote by $\Gamma^+_\al\subseteq \Gamma\setminus \Gamma^-_\al$  the set
of potential copies of $\Lambda$
that share at least one edge with $\al$, and are compatible with
$\al$. (Note that now  $\ga\in\gG^+_\ga$, unlike in \refS{se:sbc}.)

%We have already proven \eqref{eq:ctm} in \refS{Sctm};
%we next prove  (\ref{eq:ctv}). % in Theorem \ref{th:nact}, 

\begin{figure}
    \centering
    \includegraphics[scale=0.4]{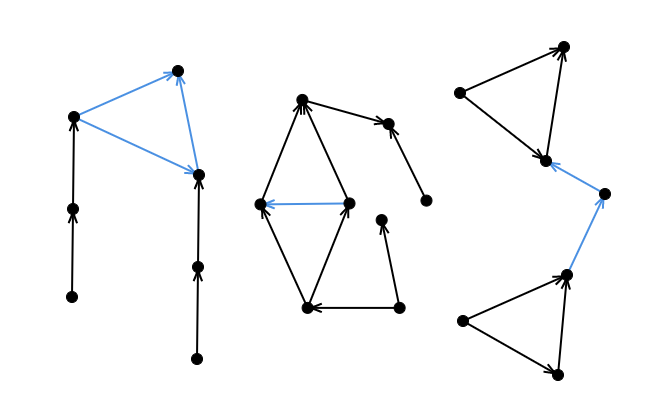}
    \caption{{\small Examples of positively correlated copies of $\Lambda$ (with $\ell=3$ and $\cT_1$ being a path with 2 edges). The shared edges are coloured blue. }}
    \label{fig:pos-balloons}
\end{figure}

%\begin{proof}[Proof of \eqref{eq:ctv}]
By edge independence, the covariance between
$\al$ and any copy of $\Lambda$ that does not belong to $\gG^+_\ga\cup\gG^-_\ga$
is zero. Thus,  
\begin{align}\label{eq:varballoon}
   \sigma_n^2&=%\sum_{\al\in\Gamma} \bigpar{\E\tone_\al - (\E\tone_\al)^2}+
 \sum_{\al\in\Gamma} \sum_{\beta\in\Gamma^+_\al} \bigpar{\E\tone_\al\tone_\beta -
          \E\tone_\al\E\tone_\beta}
 -\sum_{\al\in\Gamma} \sum_{\beta\in\Gamma^-_\al}
                \E\tone_\al\E\tone_\beta%\label{eq:varballoon0}
.\end{align}

Consider first the positive covariances, i.e., 
the first double sum in (\ref{eq:varballoon}).
We argue as in the estimate of \eqref{fq1}.
Let $\gab$ denote the union of the graphs $\ga$ and $\gb$, and note that
$\tone_\ga\tone_\gb=\tone_{\gab}$. 
Then $\gab$ is a connected graph, which either 
is unicyclic (an $\ell$-cycle with attached trees)
or multicyclic, in both cases with at most $2t$ edges outside the
2-core;
see Figure \ref{fig:pos-balloons} for examples. 
Each graph $\gab$ arises from a bounded number of pairs $(\ga,\gb)$, and thus
it follows by \refT{th:others}, 
summing over the finitely many isomorphism types of $\gab$ that can arise,
that the expected number of pairs $(\ga,\gb)$ in $G^m_n$ with
$\gb\in\gG^+_\ga$ and
$\gab$ multicyclic is
$O(\log^{2t}n)$.
Similarly, 
it follows from \eqref{eq:ctm}, 
again summing over a finite
number of
possible types of $\ga\gb$,
that the expected number of pairs 
with $\gab$ unicyclic  is $\Theta(\log^{2t+1}n)$.
Consequently, we have the exact order
\begin{align}\label{fq4}
\sum_{\al\in\Gamma} \sum_{\beta\in\Gamma^+_\al}\E\tone_\al\tone_\beta 
=\Theta(\log^{2t+1}n).   
\end{align}
Furthermore, we have 
$0\le\E\tone_\ga\E\tone_\gb\le\E\tone_\ga\tone_\gb$ for all 
$\ga$ and $\gb\in \gG^+_\ga$, with 
$\E\tone_\ga\E\tone_\gb\le\frac12\E\tone_\ga\tone_\gb$ unless
both $\ga$ and $\gb$ contain vertex 1, and it is easy to see 
that also
\begin{align}\label{fq5}
   \sum_{\al\in\Gamma} \sum_{\beta\in\Gamma^+_\al}
(\E\tone_\al\tone_\beta -\E\tone_\al\E\tone_\beta )
=\Theta(\log^{2t+1}n).
\end{align}
Alternatively, this follows from \eqref{fq4} and the argument for
$\gG^-_\ga$ below.

We see also that the dominant term is contributed 
by pairs  $\al\in\Gamma$ and  $\beta\in\Gamma^+_\al$ 
such that $\ga$ and $\gb$ have the same cycle but no other common vertices.

For the negative covariances, % in (\ref{eq:varballoon}), 
we show that the last double sum in (\ref{eq:varballoon}) is of
order $O(\log^{2t} n )$. To do so, we modify the argument for showing the sum in (\ref{fq2}) is bounded. We construct a uniform attachment
graph $G_n^{2m}$ where the first $m$ out-edges of any vertex $2\leq i\leq n$
are coloured red, and the remaining $m$ out-edges are coloured blue. 
If $\ga$ and $\gb$ are two potential copies of $\gG$,
let $\gaxb$ be the subgraph of $K^{2m}_n$ obtained by 
regarding $\ga$ as living in the red part and $\gb$ in the blue part, and
then taking the union of them 
(see Figure \ref{fig:cct} for examples). 
Note that if $\gb\in\gG^-_\ga$, then $\ga$ and $\gb$
must share at least one vertex. Thus, 
$\gaxb$ is connected; furthermore, it contains at least two cycles, and
it has at most $2t$ edges in the tree parts outside its 2-core. Once again,
with some slight abuse of notation, let $\tone_{\gaxb}$ be the
indicator that $\gaxb$ is present in $G^{2m}_n$. By independence of the edges, we have
$\E\tone_{\alpha}\E\tone_{\beta}=\E\tone_{\gaxb}$.  
Since $\Lambda$ is fixed, the total negative covariance in 
the last sum in \eqref{eq:varballoon} 
can thus be bounded by some constant times
the expected number of copies of subgraphs of the possible types of $\gaxb$ in
$G^{2m}_n$.
Again, \refT{th:others} implies that this is $O(\log^{2t}n)$.

Combining the estimates above of the  sums in \eqref{eq:varballoon}
yields the result 
$\sigma^2_n=\Theta(\log^{2t+1}n)$.
\end{proof}

\begin{proof}[Proof of \refT{th:others}]
We have proved \eqref{cq0} in \refS{Sctm},
and \eqref{cq01}
follows from \refT{th:bdexpsub} and an argument 
similar to the proof above of \eqref{eq:ctv}
(where we may simplify and ignore the negative terms),
    noting that a subgraph obtained by merging two multicyclic and leaf-free
    subgraphs at one or more vertices and edges is still multicyclic and
    leaf-free. 
\end{proof}

\subsection{Proof of \refT{th:nact}}

We now complete the proof of \refT{th:nact} by 
proving the normal approximation result \eqref{eq:ctna}.

\begin{proof}[Proof of \eqref{eq:ctna}]
We bound the error terms appearing in (\ref{eq:wassersteinbd}) using the
coupling described in Section \ref{se:sbc}. 
Let $W^s_n$ be the size-bias
version of $W_n$ defined there.
For $\al\in\gG$, 
let
$\Gamma^-_\al$ and
$\Gamma^+_\al$ be the subsets of $\gG$ defined at the beginning of this
subsection,
and
let as in \refS{se:sbc}, $G^{\al}$ be the graph $G$
forced to contain all edges of $\ga$.
Denote by $\tone_\beta^\al$ the indicator
that a copy $\beta$ of $\gL$ is in $G^{\al}$, and let 
$\ttone^\ga_\gb:=\tone^\ga_\gb-\tone_\gb$.
Also,
if $\gb\in\gG^+_\ga$, let $\gb\setminus\ga$ be the graph obtained from
$\gb$ by deleting the edges that also are in $\ga$; then
$\tone^\ga_\gb=\tone_{\gb\setminus\ga}$. 

If $\beta \in \Gamma^+_\al$, then $\ttone^\ga_\gb\ge0$ 
since if $\gb$ exists in $G$, then it will also exist in $G^\ga$;
furthermore,
if $\beta \in \Gamma^-_\al$, then 
$\tone^\ga_\gb=0$ and thus $\ttone^\ga_\gb=-\tone_\gb$, and if 
$\gb\notin\gG^+_\ga\cup\gG^-_\ga$, then $\ttone^\ga_\gb=0$. 
Hence, the construction of $W^s$ in \refS{se:sbc} yields
    \begin{align*}
        \E[W^s_n-W_n\mid G]&
= \sum_{\al\in \Gamma}\frac{\E\tone_\al}{\mu_n}\sum_{\gb\in\gG}\ttone^\ga_\gb
= \sum_{\al\in \Gamma}\frac{\E\tone_\al}{\mu_n}
\biggpar{\sum_{\beta\in \Gamma^+_\al} \ttone^\al_{\beta} - \sum_{\beta\in\Gamma^-_\al}  
\tone_\beta }
\\&
=: \sum_{\al\in \Gamma}\frac{\E\tone_\al}{\mu_n}
\bigpar{R^+_\ga-R^-_\ga}.
\numberthis\label{qa0}
    \end{align*}
Taking the variance yields, using the (Cauchy--Schwarz) inequality $(a+b)^2\le
2(a^2+b^2)$ for all real $a$ and $b$,
\begin{align}\label{qa1}
 &   \mu_n^2\var \bigpar{\E[W^s_n-W_n\mid G]}
%\\&
\le 2\Var\biggpar{\sum_{\al\in \Gamma}\E\tone_\al R^+_\ga}
+2\Var\biggpar{\sum_{\al\in \Gamma}\E\tone_\al R^-_\ga}.
\end{align}
We proceed to bound the variances on the right-hand side.
The first equals
\begin{align*}\numberthis\label{qa2}
   &
A_1:=\sum_{\alpha\in\Gamma} \sum_{\beta\in\Gamma}\E\tone_\al\E\tone_\beta
     \cov(R^+_\al,R^+_\beta)  
%\\ &
=\sum_{\alpha\in\Gamma} \sum_{\beta\in\Gamma} 
\sum_{\gamma_1\in\Gamma^+_\al}\sum_{\gamma_2\in\Gamma^{+}_\beta} 
\E\tone_\al\E\tone_\beta
\Cov\bigpar{\ttone^\al_{\gamma_1},\ttone^\beta_{\gamma_2}}
.\end{align*}
In \eqref{qa2}, we only need to consider pairs of $\gamma_1,\gamma_2$ that
have at least one edge in common, 
since otherwise $\ttone^\al_{\gamma_1}$ and $\ttone^\beta_{\gamma_2}$ are
independent. %, and thus their covariance is zero. 
Moreover, since $\gamma_1\in\Gamma^+_\al$ and $\gamma_2\in\Gamma^{+}_\beta$,
$\gamma_1$ has an edge (and thus a vertex) in common with $\ga$
and
$\gamma_2$ has an edge (and thus a vertex) in common with $\gb$.
We argue similarly to the proof of \eqref{eq:ctv} above, this time
considering the uniform attachment graph $G^{3m}_n$, which we regard as three
independent copies of $G$ that are coloured red, blue, and green.
Let $\gabxcc$ be the graph formed by a red copy of $\ga$, a blue copy of
$\gb$, and a green copy of 
$(\gamma_1\setminus\ga)\cup(\gamma_2\setminus\gb)$.
Then
\begin{align}\label{qa3}
\Cov\bigpar{\ttone^\al_{\gamma_1},\ttone^\beta_{\gamma_2}}
\le \E\bigsqpar{\ttone^\al_{\gamma_1}\ttone^\beta_{\gamma_2}}
\le \E\bigsqpar{\tone^\al_{\gamma_1}\tone^\beta_{\gamma_2}}
= \E\bigsqpar{\tone_{\gamma_1\setminus\ga}\tone_{\gamma_2\setminus\gb}}
\end{align}
and thus
\begin{align}\label{qa4}
\E\tone_\ga\E\tone_\gb\Cov\bigpar{\ttone^\al_{\gamma_1},\ttone^\beta_{\gamma_2}}
\le 
\E\tone_\ga\E\tone_\gb
\E\bigsqpar{\tone_{\gamma_1\setminus\ga}\tone_{\gamma_2\setminus\gb}}
=
\E\bigsqpar{\tone_\gabxcc}
.\end{align}
Note that $\gabxcc$ is a connected graph, which has at least one cycle. 
(It unicyclic when the graphs $\ga,\gb,\gamma_1,\gamma_2$ have the same cycle,
and their attached trees are disjoint.) 
The edges in the trees attached to the 2-core of $\gabxcc$ come  from the
attached trees in $\ga,\gb,\gamma_1$, and $\gamma_2$, and thus the total
number of them is at most $4t$. 
Consequently, 
\eqref{eq:ctm} and \refT{th:others} apply
and show together with \eqref{qa2}--\eqref{qa4},
summing
over the finite number of possible types of the graph $\gabxcc$,
\begin{align}\label{qa5}
  A_1= O\bigpar{\log^{4t+1}n}.
\end{align}

Similarly,
the second variance on the right-hand side of \eqref{qa1} equals
\begin{align*}\numberthis\label{qb2}
A_2:=\sum_{\alpha\in\Gamma} \sum_{\beta\in\Gamma}\E\tone_\al\E\tone_\beta
     \cov(R^-_\al,R^-_\beta)  
%\\ &
=\sum_{\alpha\in\Gamma} \sum_{\beta\in\Gamma} 
\sum_{\gamma_1\in\Gamma^-_\al}\sum_{\gamma_2\in\Gamma^{-}_\beta} 
\E\tone_\al\E\tone_\beta
\Cov\bigpar{\tone_{\gamma_1},\tone_{\gamma_2}}
.\end{align*}
Again, we only need to consider pairs of $\gamma_1,\gamma_2$ that
have at least one edge in common.
%since otherwise $\ttone^\al_{\gamma_1}$ and $\ttone^\beta_{\gamma_2}$ are
%independent. %, and thus their covariance is zero. 
Moreover, since $\gamma_1\in\Gamma^-_\al$ and $\gamma_2\in\Gamma^{-}_\beta$,
$\gamma_1$ has a vertex in common with $\ga$
and
$\gamma_2$ has a vertex in common with $\gb$.
We consider again $G^{3m}_n$ and now let $\gabcc$ be the graph consisting of
a red copy of $\ga$, a blue copy of $\gb$, and green copies of $\gam_1$ and
$\gam_2$. 
Then
\begin{align}\label{qb4}
\E\tone_\ga\E\tone_\gb\Cov\bigpar{\tone_{\gamma_1},\tone_{\gamma_2}}
\le 
\E\tone_\ga\E\tone_\gb
\E\bigsqpar{\tone_{\gamma_1}\tone_{\gamma_2}}
=
\E\bigsqpar{\tone_\gabcc}
,\end{align}
and $\gabcc$ is a connected graph with at least 
three cycles
and at most $4t$
edges outside its 2-core.
Consequently, we now obtain from
\eqref{eq:ctm} and \refT{th:others}
\begin{align}\label{qb5}
  A_2= O\bigpar{\log^{4t}n}.
\end{align}

By \eqref{qa1} and the bounds \eqref{qa5} and \eqref{qb5} above, 
we deduce that 
\begin{align}\label{qc1}
\mu_n^2\var\bigpar{\E[W^s_n-W_n\mid G]}=O\bigpar{\log^{4t+1}n}.   
\end{align}
Moreover, $W_n$ is determined by $G$, and thus the conditional
Jensen's inequality yields
\begin{align}\label{qc2}
  \Var\bigpar{\E[W^s_n-W_n\mid W_n]}
=
  \Var\bigpar{\E\bigsqpar{\E[W^s_n-W_n\mid G]\mid W_n}}
\le
  \Var\bigpar{\E[W^s_n-W_n\mid G]}.
\end{align}
Noting also that $\sigma^2_n=\Theta(\logn{2t+1})$ by \eqref{eq:ctv}, 
the first error term in (\ref{eq:wassersteinbd}) is therefore
\begin{align}\label{qc3}
\frac{\mu_n}{\sigma^2_n}\sqrt{\frac{2}{\pi}}\sqrt{\var\bclr{\E[W^s_n-W_n\mid W_n]}}
=O\bigpar{\logn{(4t+1)/2-(2t+1)}}
=O\bigpar{\log^{-\frac{1}{2}} n}.
\end{align}

We now turn to the second error term in (\ref{eq:wassersteinbd}). 
By first conditioning on the graph $G$, we have in analogy with \eqref{qa0},
    \begin{align*}\numberthis\label{qd0}
        \E\bigsqpar{(W^s_n-W_n)^2\mid G}&
= \sum_{\al\in \Gamma}\frac{\E\tone_\al}{\mu_n}
\bigpar{R^+_\ga-R^-_\ga}^2
\le2 \sum_{\al\in \Gamma}\frac{\E\tone_\al}{\mu_n}
\bigpar{(R^+_\ga)^2+\xpar{R^-_\ga}^2}
,\end{align*}
and thus
    \begin{align*}\numberthis\label{qd1}
\mu_n\E\bigsqpar{(W^s_n-W_n)^2}&
\le2 \sum_{\al\in \Gamma}\E\tone_\al(R^+_\ga)^2
+2 \sum_{\al\in \Gamma}\E\tone_\al(R^-_\ga)^2
.\end{align*}
We argue as above.
The first sum on the right-hand side equals
\begin{align*}\numberthis\label{qd2}
B_1:=\sum_{\alpha\in\Gamma} \E\tone_\al
\sum_{\gamma_1\in\Gamma^+_\al}\sum_{\gamma_2\in\Gamma^{+}_\ga} 
\Cov\bigpar{\ttone^\al_{\gamma_1},\ttone^\ga_{\gamma_2}}
.\end{align*}
%Again, it suffices to consider pairs of $\gamma_1,\gamma_2$ that
%have at least one edge in common.
Hence, using \eqref{qa3},
\begin{align*}
B_1
\le
\sum_{\alpha\in\Gamma} \sum_{\gamma_1\in\Gamma^+_\al}\sum_{\gamma_2\in\Gamma^{+}_\ga} 
\E\tone_\al
\E\bigsqpar{\tone_{\gamma_1\setminus\ga}\tone_{\gamma_2\setminus\ga}}
=
\sum_{\alpha\in\Gamma} \sum_{\gamma_1\in\Gamma^+_\al}\sum_{\gamma_2\in\Gamma^{+}_\ga} 
\E\tone_\gaxcc
\numberthis\label{qd3}
,\end{align*}
where 
$\gaxcc$ is  the subgraph formed by a red copy of $\ga$
and a blue copy of 
$(\gamma_1\setminus\ga)\cup(\gamma_2\setminus\ga)$,
regarded as a subgraph of $K^{2m}_n$  coloured red and blue as above,
and $\tone_\gaxcc$ is the indicator that this subgraph exists in $G^{2m}_n$.
Since $\gamma_1,\gamma_2\in\gG_\ga^+$, they have at least one vertex each in
common with $\ga$, and thus $\gaxcc$ is connected. Moreover, $\gaxcc$ has at
least one cycle, and at most $3t$ edges outside its 2-core.
Consequently, we now obtain from
\eqref{eq:ctm} and \refT{th:others}
\begin{align}\label{qd5}
  B_1= O\bigpar{\log^{3t+1}n}.
\end{align}

Similarly,
the second sum in \eqref{qd1} equals
\begin{align}\label{qe2}
B_2&:=\sum_{\alpha\in\Gamma} \E\tone_\al
\sum_{\gamma_1\in\Gamma^-_\al}\sum_{\gamma_2\in\Gamma^{-}_\ga} 
\Cov\bigpar{\tone_{\gamma_1},\tone_{\gamma_2}}
%.\end{align*}
%\begin{align*}
%B_2
\le
\sum_{\alpha\in\Gamma} \sum_{\gamma_1\in\Gamma^-_\al}\sum_{\gamma_2\in\Gamma^{-}_\ga} 
\E\tone_\al
\E\bigsqpar{\tone_{\gamma_1}\tone_{\gamma_2}}
\notag\\&\phantom:
=
\sum_{\alpha\in\Gamma} \sum_{\gamma_1\in\Gamma^-_\al}\sum_{\gamma_2\in\Gamma^{-}_\ga} 
\E\tone_\gacc
,\end{align}
where 
$\gacc$ is  the subgraph formed by a red copy of $\ga$
and a blue copy of 
$\gamma_1\cup\gamma_2$.
%regarded as a subgraph of $K^{2m}_n$  coloured red and blue as above,
%and $\tone_\gaxcc$ is the indicator that this subgraph exists in $G^{2m}_n$.
%Since $\gamma_1,\gamma_2\in\gG_\ga^-$, they have at least one vertex each in
%common with $\ga$, and thus $\gaxcc$ is connected. Moreover, $\gaxcc$ has at
%least one cycle, and at most $3t$ edges outside its 2-core.
We  obtain from
\eqref{eq:ctm} and \refT{th:others}
\begin{align}\label{qe5}
  B_2= O\bigpar{\log^{3t+1}n}.
\end{align}

We obtain from \eqref{qd1}, \eqref{qd5}, and \eqref{qe5} that
\begin{align}\label{qf1}
\mu_n\E\bigsqpar{(W^s_n-W_n)^2}&
= O\bigpar{\log^{3t+1}n}.
\end{align}
By this and \eqref{eq:ctv}, we conclude that
the second error term in (\ref{eq:wassersteinbd}) is 
\begin{align}\label{qf3}
     \frac{\mu_n}{\sigma^3_n}\E[(W^s_n-W_n)^2]
=O\bigpar{\logn{3t+1-3(2t+1)/2}}
=O\bigpar{\log^{-\frac{1}{2}} n}.
\end{align}

Finally, we use \eqref{qc3} and \eqref{qf3} in \eqref{eq:wassersteinbd}
and obtain  \eqref{eq:ctna}, 
which completes the proof of \refT{th:nact}.
\end{proof}

%%XXXXXXXXXXXXXXXXXXXXXXXXXXXXXXXXXXXXXXXXXXXXXXXXXXXXXXXXXXX
%%\input{misc}
\section{Proof of the normal approximations for  trees}\label{se:natrees}
In this section we prove Theorem \ref{th:natree} and  \refT{th:na2star}. 

\subsection{Proof of Theorem \ref{th:natree}.} 
We start by proving \eqref{eq:treem}:
By Lemma \ref{le:gentree}(i) and \eqref{pga}, an upper bound is,
by summing in the order $k_1,\dots,k_t$,
\begin{align}
  \mu_n
\le 
C \sum_{1\le k_1<\dots<k_t\le n}\prod^t_{i=2}\frac{1}{k_i-1}
\le C \sum_{k_t\le n}1=
Cn,
\end{align}
The lower bound follows similarly, or simply because there are
$\Theta(n^t)$ potential copies of $\gL$, and each has probability $\ge
n^{-(t-1)}$. 
    
\eqref{eq:treev}:
To obtain an upper bound on $\sigma^2_n$, we argue as in the proof of
\eqref{eq:ctv}. We use again
\eqref{eq:varballoon}, and the only difference from the argument
in \refS{se:naballoons} is that 
for $\al\in\gG$ and
$\beta\in \Gamma^{+}_\ga$, their union $\al\beta$ 
does not have to contain
a cycle; however, it is always
a connected graph with at most $2t$ edges. 
There is still a finite number of types of $\ga\gb$, 
and we obtain by using (depending on the number of cycles in $\ga\gb$)
\eqref{eq:treem},
\eqref{eq:ctm} and \refT{th:others} 
\begin{align}
\sum_{\ga\in\gG}\sum_{\gb\in\gG^+_\ga}\E\tone_{\ga\gb}
\le O(n)+O\bigpar{\log^{2t}n}
= O(n),
\end{align}
which implies \eqref{eq:treev}.

If $\sigma^2_n=\Theta(n)$, then we can use Theorem \ref{th:wassersteinbd} to
prove the asymptotic normality of $W_n$ in a similar vein as before. 
Again, we argue as in the proof of \refT{th:nact} in \refS{se:naballoons},
and the only difference is that the graphs $\gabxcc$, $\gabcc$, $\gaxcc$,
$\gacc$ do not have to contain a cycle. As in the estimation of variance
above, we therefore use not only
\eqref{eq:ctm} and \refT{th:others} but also \eqref{eq:treem} in the estimates;
as a result we obtain the bound $O(n)$ in 
%\eqref{qa5}, \eqref{qb5},
\eqref{qc1}
%, \eqref{qd1}, \eqref{qe5} 
and \eqref{qf1}.
Hence,
if $\sigma_n^2=\Theta(n)$, the terms in (\ref{eq:wassersteinbd}) are 
    \begin{equation}
        \frac{\mu_n}{\sigma_n^2}\sqrt{\var(\E[W^s_n-W_n\mid G])} = O(n^{-1/2}),\quad  \frac{\mu_n}{\sigma_n^3} \E[(W^s_n-W_n)^2]=O(n^{-1/2}),
    \end{equation}
 implying (\ref{eq:treena}). \qed

\subsection{Proof of \refT{th:na2star}.}\label{SP2} 
Recall that the  upper bound $O(n)$ 
for the variance $\gss_n$ 
is already proved in \eqref{eq:treev}.
To obtain a lower bound, 
we state and
prove the following general lemma. 
In our application of it,
$X$ and $\xi_\nu$ will be taken as the count of
trees $W_n$ and the recipient of a directed edge $\nu$ in the uniform
attachment graph $G^m_n$.

\begin{lemma}\label{le:varlb}
Let $(\xi_\nu)_{\nu\in\cI}$ be a family of independent random variables
(where $\cI$ is an arbitrary index set), and
let $X$ be any random variable such that $\E X^2<\infty$.
For each $\nu\in\cI$, let
\begin{align}\label{a1}
  X_\nu:=\E(X\mid\xi_\nu)-\E X.
\end{align}
Then
\begin{align}\label{a2}
  \var (X) \ge \sum_{\nu\in\cI} \var (X_\nu)
= \sum_{\nu\in\cI} \var \big[\E (X\mid \xi_\nu)\big].
\end{align}
\end{lemma}

\begin{proof}
Assume first that $\cI$ is finite, and define
\begin{align}\label{a3}
  Y:=\suma X_\nu
\end{align}
Note that 
\begin{align}\label{a3.5}
 \E Y=\suma \E X_\nu=0.
\end{align}

It follows from \eqref{a1} that the random variable $X_\nu$ is a function of
$\xi_\nu$. 
Consequently, the variables $X_\nu$ are independent, and thus
\begin{align}\label{a4}
\E (Y^2)=
  \var(Y)=\suma\var(X_\nu).
\end{align}
Moreover, \eqref{a1} implies also
\begin{align}\label{a5}
  \E\bbcls{(X-\E X)X_\nu}
=   \E\bbcls{\E(X-\E X\mid\xi_\nu)X_\nu}
=\E\bclr{X_\nu^2}=\var(X_\nu).
\end{align}
Consequently,
\begin{align}\label{a6}
\E\bbcls{(X-\E X)Y}
=\suma  \E\bbcls{(X-\E X)X_\nu}
=\suma\var(X_\nu)
=\E (Y^2).
\end{align}
Hence,
\begin{align}\label{a7}
\E\bbcls{(X-\E X-Y)Y}
=0
\end{align}
and thus
\begin{align}\label{a8}
\E\bbcls{(X-\E X)^2}&
=\E\bbcls{(X-\E X-Y)^2}
+2\E\bbcls{(X-\E X-Y)Y}
+\E\bbcls{Y^2}
\notag\\&
=\E\bbcls{(X-\E X-Y)^2}
+\E\bbcls{Y^2}
\notag\\&
\ge
\E\bcls{Y^2}.
\end{align}
The result follows from \eqref{a8} and \eqref{a4}, which completes the proof
for finite $\cI$. 

If $\cI$ is infinite, the case just proved shows that for every finite
subset $\cI_1\subset\cI$, we have
\begin{align}\label{a9}
  \var(X)\ge\sum_{\nu\in\cI_1}\var(X_\nu),
\end{align}
which implies \eqref{a2}.
\end{proof}
\begin{remark}
From a more abstract point of view, 
a minor modification of the proof above shows that
$Y$ is the orthogonal projection of $X-\E X$ onto
the linear  subspace of $L^2$ consisting of sums of the type 
$\suma f_\nu(\xi_\nu)$, and the inequality \eqref{a2} is thus
an instance of the fact that orthogonal projections in $L^2$ never increase
the variance (or the norm).
\end{remark}

\begin{proof}[Proof of \refT{th:na2star}]
Fix $j\in[n]$ and $a\in[m]$. 
Recall that the edge $j\aaa$ in $G^m_n$ is randomly chosen as one of the
edges $ji\aaa$ ($i\in[j-1]$) in $K^m_n$; we denote the recipient of
the edge $j\aaa$ by $\xi=\xi_{j,a}$; thus 
$\xi$ is uniformly distributed on $[j-1]$, and
$j\aaa=j\xi\aaa$.
We condition on $\xi$, and decompose $W_n$ as
\begin{align}
  W_n=\WWn0+\WWn1+\WWn2,
\end{align}
where 
\begin{itemize}
\item 
$\WWn0$ is the number of copies of $\Sl$ that do not contain the
edge $ji\aaa$ for any $i\in[j-1]$.
\item 
$\WWn1$ is the number of copies $\ga$ of $\Sl$ that contain
$j\aaa=j\xi\aaa$ and such that the center of $\ga$ is $j$.
\item 
$\WWn2$ is the number of copies $\ga$ of $\Sl$ that contain
$j\aaa=j\xi\aaa$ and such that the center of $\ga$ is $\xi$.
\end{itemize}
Then $\WWn0$ is clearly independent of $\xi$.
Moreover, for every $i<j$,
$\E\xpar{\WWn1\mid\xi=i}$ is the expected number of stars $\Sli$ with center $j$
and $\ell-1$ leaves in $[n]\setminus\set{i,j}$.
Since the edges from $j$ have their endpoints uniformly distributed in
$[j-1]$, this number does not depend on $i$.
In other words, $\E(\WWn1\mid\xi)$ does not depend on $\xi$.

Similarly, 
$\E\xpar{\WWn2\mid\xi=i}$ is the expected number of stars $\Sli$ with center $i$
and $\ell-1$ leaves in $[n]\setminus\set{i,j}$.
Hence,
\begin{align}
  \E\bigpar{\WWn2\mid\xi=i}=\sum_{q=0}^{\ell-1}\WIJ q,
\end{align}
where $\WIJ q$ is the number of such copies of $\Sli$ with $q$ leaves in 
$(i,n]$ and $r:=\ell-1-q$ in $[1,i)$.
Assume for simplicity that $i\ge n/10$.
Then, counting the number of ways to choose first the vertices and then the
edges of such a copy, and multiplying with the probability that it exists in
$G^m_n$, 
\begin{align}
  \E\WIJ q &
= \binom{i-1}r\sum_{\substack{i<i_1<\dots<i_q\le n\\i_1,\dots,i_q\neq j}} 
\binom{m}{r} r!\, m^q(i-1)^{-r}\prod_{k=1}^q(i_k-1)^{-1}
\notag\\&
=\Bigpar{1+O\Bigpar{\frac{1}n}} 
\binom{m}{r}  m^q\frac{1}{q!} 
\lrpar{\sum_{p=i+1}^n\frac{1}{p-1}+O\Bigpar{\frac{1}{n}}}^q
\notag\\&
=
\binom{m}{r}  \frac{m^q}{q!} \log^q \Bigpar{\frac{n}{i}}
+O\Bigpar{\frac{1}n}
.\end{align}
Since $\log(n/i)$ is positive and monotonically decreasing in
$i\in[n/10,n)$,
it follows that if $\frac{n}{10}\le i\le k< j\le n$, then
\begin{align}
  \E\bigpar{W_n\mid\xi=i} - \E\bigpar{W_n\mid\xi=k}&
=  \E\bigpar{\WWn2\mid\xi=i} - \E\bigpar{\WWn2\mid\xi=k}
\notag\\&
=\sum_{q=0}^{\ell-1}\bigpar{\WIJ{q}-\WKJ{q}}
\notag\\&
=\sum_{q=0}^{\ell-1}\binom{m}{r}\frac{m^q}{q!} 
\Bigpar{\log^q \Bigpar{\frac{n}{i}}-\log^q \Bigpar{\frac{n}{k}}
+O\Bigpar{\frac{1}n}}
\notag\\&
\ge \frac{m^{\ell-1}}{(\ell-1)!} 
\Bigpar{\log^{\ell-1} \Bigpar{\frac{n}{i}}-\log^{\ell-1}\Bigpar{ \frac{n}{k}}}
+O\Bigpar{\frac{1}n}
.\end{align}
In particular, if
$i\in[n/10,n/4]$ and $k\in[n/3,n/2]$, %and $j\in(n/2,n]$,
then
\begin{align}\label{GIIA}
  \E\bigpar{W_n\mid\xi=i} - \E\bigpar{W_n\mid\xi=k}&
\ge \frac{m^{\ell-1}}{(\ell-1)!} \Bigpar{\log^{\ell-1} 4-\log^{\ell-1}3}
+O\Bigpar{\frac{1}n}
\ge c
\end{align}
for some $c>0$ if $n$ is large enough,
say $n\ge n_1$.
 
Let $\xi'$ be an independent copy of $\xi=\xi_{j,a}$ and write
$g(\xi):=\E(W_n\mid\xi)$.
Assume $n\ge n_1$.
Then \eqref{GIIA} shows that $g(i)-g(k)\ge c$ if 
$i\in[n/10,n/4]$ and $k\in[n/3,n/2]$;
hence it follows that if
$j>n/2$, then
\begin{align}\label{eq:qh1}
  \Var\bigsqpar{ \E(W_n\mid\xi)}
&=\Var g(\xi)
=\tfrac12\E\lrpar{g(\xi)-g(\xi')}^2
\notag\\&
\ge \tfrac12\P\bigpar{\xi\in[n/10,n/4]}\P\bigpar{\xi'\in[n/3,n/2]} c^2
\ge c_1,
\end{align}
for some constant $c_1>0$.

We now apply Lemma \ref{le:varlb} with the family of all random variables
$\xi_{j,a}$, thus letting $\cI=\set{2,\dots,n}\times[m]$.  
Then \eqref{a2} and the lower bound in \eqref{eq:qh1}
(for $j\ge n/2$)
show that  for $n\geq n_1$,
\begin{align}
    \var(W_n) \geq \sum^n_{j=\ceil{n/2}} \sum^m_{a=1} c_1 = \Omega(n).
\end{align}
Since $\Var W_n = \sigma_n^2=O(n)$ by (\ref{eq:treev}), this completes the proof.
\end{proof}

%%XXXXXXXXXXXXXXXXXXXXXXXXXXXXXXXXXXXXXXXXXXXXXXXXXXXXXXXXXXX
%\input{other-non-tree-graphs}

\section{Proof of \refT{Ttail}}

\begin{proof}[Proof of \refT{Ttail}]
Let $h'$ be the number of vertices in $H'$, and let $r\in[h]$ be the number
of them such that attaching an edge at that vertex yields a copy of $H$.

Fix a copy $H'_0$ of $H'$ in $G^m_\infty$, and let its vertices
be $k_1<\cdots<k_{h'}$.
For every $j\in\bbN$ and $a\in[m]$, let $I_{ja}$ be the indicator of the
event that the edge $j^{(a)}$ in $G^n_m$ has exactly one endpoint in $H'_0$ 
and moreover that attaching $j^{(a)}$ to $H'_0$ yields a copy of $H$.
Then, for $n\ge k_{h'}$,  the number of copies of $H$ in $G^m_n$ that contain
the given subgraph $H'_0$ is 
\begin{align}
W_n^{H'_0}:=\sum_{j=1}^n\sum_{a=1}^m I_{ja}  
=
\sum_{j = k_{h'}+1}^n\sum_{a=1}^m I_{ja}  +O(1).
\end{align}
Condition 
on the existence of $H_0'$.
Then the variables $I_{ja}$ in the final sum are independent and have the 
Bernoulli distributions $\Be(r/(j-1))$.
It follows by the lemma below  that 
as \ntoo,
$W_n^{H'_0}/\log n \asto rm$.
Thus, every copy of $H'$ in $G^m_\infty$ has asymptotically $rm\log n$
attached edges that make copies of $H$.
Summing over the a.s.\ finite number of copies $H'_0$ of $H'$ in $G^m_\infty$,
we obtain
\begin{align}
  \frac{W_n}{\log n}\asto rm W'_\infty,
\end{align}
which verifies \refConj{conj:A-B} in this case.
\end{proof}

The proof used the following simple law of large numbers; it is certainly
known but for completeness we include a proof.
\begin{lemma}\label{LLN}
Let $I_i\in\Be(p_i)$ be independent random variables with $p_i\in[0,1]$,
for $i=1,2,\dots$.
If\/ $\sum_{i=1}^\infty p_i=\infty$, then
\begin{align}\label{lln}
  \frac{\sumin I_i}{\sumin p_i}\asto 1.
\end{align}
\end{lemma}
\begin{proof}
Let $b_n:=\sumin p_i$.
  We have
  \begin{align}
    \sumn \frac{\E(I_n-p_n)^2}{b^2_n}&
\le \sumn \frac{p_n}{b_n^2}
=\frac{1}{p_1}+\sum_{n=2}^\infty \frac{b_n-b_{n-1}}{b_n^2}
\le \frac{1}{p_1}+\sum_{n=2}^\infty \int_{b_{n-1}}^{b_n}\frac{\ddx s}{s^2}
\notag\\&
=\frac{1}{p_1}+ \int_{b_{1}}^\infty\frac{\ddx s}{s^2}
<\infty.
  \end{align}
Since the $I_n$ are independent, it follows that
$b_n^{-1}\sumin (I_i-p_i)\asto0$, 
see \cite[Theorem VII.8.3]{FellerII},
which is equivalent to \eqref{lln}.
\end{proof}

\section*{Acknowledgement}
Tiffany Y.\ Y.\ Lo thanks
Nathan Ross for pointing her to some results in
Stein's method.

%%XXXXXXXXXXXXXXXXXXXXXXXXXXXXXXXXXXXXXXXXXXXXXXXXXXXXXXXXXXX
%\input{other-non-tree-graphs}

%\nocite{*}
\printbibliography

\end{document}